\documentclass[11pt]{amsart}
\usepackage{amsmath,amssymb,amsthm,mathtools}
\usepackage[margin=1.15in]{geometry}
\usepackage{enumitem}
\usepackage[colorlinks=true,linkcolor=blue,citecolor=blue,urlcolor=blue]{hyperref}

\newtheorem{theorem}{Theorem}[section]
\newtheorem{proposition}[theorem]{Proposition}
\newtheorem{lemma}[theorem]{Lemma}
\newtheorem{corollary}[theorem]{Corollary}
\theoremstyle{definition}
\newtheorem{assumption}[theorem]{Assumption}
\theoremstyle{remark}
\newtheorem{remark}[theorem]{Remark}
\numberwithin{equation}{section}

\newcommand{\R}{\mathbb{R}}
\newcommand{\cM}{\mathcal{M}}
\newcommand{\cS}{\mathcal{S}}
\newcommand{\cL}{\mathcal{L}}
\newcommand{\cF}{\mathcal{F}}
\newcommand{\cH}{\mathcal{H}}
\newcommand{\vp}{\varpi}
\DeclareMathOperator{\supp}{supp}

\begin{document}

\title[Non-local operators with square-Dini coefficients]{Boundedness and $L_p$-estimates for non-local operators with square-Dini continuous coefficients}

\author[H. Dong]{Hongjie Dong}
\address[H. Dong]{Division of Applied Mathematics, Brown University, 182 George Street, Providence, RI 02912, USA}
\email{Hongjie\_Dong@brown.edu}
\thanks{H. Dong was partially supported by the NSF under agreement DMS-2350129.}

\author[P. Jung]{Pilgyu Jung}
\address[P. Jung]{Department of Mathematics, Yonsei University, 50 Yonsei-Ro, Seodaemun-Gu, Seoul, 03722, Republic of Korea}
\email{pilgyujung@yonsei.ac.kr}
\thanks{P.  Jung was supported by the Basic Science Research Program through the National Research Foundation of Korea (NRF) funded by the Ministry of Education (RS-2026-25572931) and was partially supported by an NRF grant funded by the Korean government (MSIT) (RS-2024-00408003).}

\author[D. Kim]{Doyoon Kim}
\address[D. Kim]{Department of Mathematics, Korea University, 145 Anam-ro, Seongbuk-gu, Seoul, 02841, Republic of Korea}
\email{doyoon\_kim@korea.ac.kr}
\thanks{D. Kim was supported by the National Research Foundation of Korea (NRF) grant funded by the Korea government (MSIT) (RS-2025-16065192) and by a Korea University Grant.}

\subjclass[2020]{35R11, 35B45, 47B38, 42B25}

\keywords{non-local equations, $x$-dependent kernels, square-Dini condition, boundedness of operators, $L_p$-estimates, Littlewood--Paley decomposition}

\begin{abstract}
We prove the boundedness of the non-local operator $\cL^a$ of order $\sigma\in(0,2)$ with kernel $a(x,y)|y|^{-d-\sigma}$ from $H^\sigma_{p,w}(\R^d)$ to $L_{p,w}(\R^d)$ for the whole range of $p\in(1,\infty)$, where $w$ is a Muckenhoupt weight. The coefficient $a(x,y)$ is bounded, merely measurable in $y$, and its modulus of continuity $\omega$ in $x$ satisfies the square-Dini condition $\int_0^1\omega(t)^2t^{-1}\,dt<\infty$. This condition allows, for instance, $\omega(t)=(\log(1/t))^{-1}$, which is not a Dini modulus, and it improves the H\"older continuity assumption in previous work. Using the boundedness result, we obtain global estimates and the unique solvability in weighted $L_p$ spaces, as well as interior $L_p$-estimates, for the corresponding elliptic and parabolic equations. We also show that the square-Dini condition is sharp. By adapting a recent construction of C.~Liu, we show that for every regular modulus $\omega$ with $\int_0^1\omega(t)^2t^{-1}\,dt=\infty$, there is a coefficient with modulus of continuity $\omega$ for which the boundedness, the global estimates, and the interior estimates all fail. In particular, on the scale $\omega(t)=(\log(1/t))^{-\theta}$, the threshold is $\theta=1/2$.
\end{abstract}

\maketitle

\section{Introduction}

We consider non-local elliptic equations of the form
\begin{equation}\label{eq:ell}
\lambda u-\cL^a u=f\quad\text{in }\R^d,
\end{equation}
and the corresponding parabolic equations $u_t-\cL^au+\lambda u=f$, where the non-local operator $\cL^a$ is defined by
\begin{equation}\label{eq:L}
\cL^a u(x)=\int_{\R^d}\big(u(x+y)-u(x)-\chi_\sigma(y)(\nabla u(x),y)\big)\,a(x,y)\,\frac{dy}{|y|^{d+\sigma}},
\end{equation}
$\sigma\in(0,2)$, and $\chi_\sigma=1_{\sigma>1}+1_{\sigma=1}1_{B_1}$. The coefficient $a(x,y)$ is bounded and measurable. Non-local operators as above appear in the equations describing various phenomena in physics, mathematical finance, biology, and fluid dynamics; see, for instance, \cite{CV10,CT04}. In this paper, we are concerned with the following question: 
\begin{quote}
\itshape
What regularity of $a$ in $x$ is needed for the boundedness
of $\cL^a$ from $H^\sigma_p(\R^d)$ to $L_p(\R^d)$
and for the $L_p$-theory of \eqref{eq:ell}?
\end{quote}
For local operators, the situation is well understood. If $\cL u=a^{ij}D_{ij}u$, the boundedness of $\cL$ from $W^2_p$ to $L_p$ requires no regularity of the coefficients, and the $W^2_p$-estimates hold when the coefficients are merely continuous, or even of vanishing mean oscillation; see, for instance, \cite{GT01,Kry08}. For non-local operators, the boundedness is far from obvious even when $a=a(y)$. In early work, some conditions were imposed on $a(y)$, such as symmetry or homogeneity and smoothness in $y$; see \cite{BB07,MP92}. These restrictions were removed in \cite{DK12}, where $a(y)$ is merely bounded and measurable. The boundedness of non-local operators with $x$-dependent coefficients was first considered in \cite{MP14} for sufficiently large $p$, when $a$ is H\"older continuous in $x$. In \cite{DJK23}, the restriction on $p$ was removed and weighted spaces with Muckenhoupt weights were considered, again under the assumption that $a$ is H\"older continuous in $x$ with an arbitrarily small exponent. See also \cite{Zha13} for the $L_p$-maximal regularity of non-local parabolic equations. The question of whether the boundedness and the $L_p$-estimates still hold when $a$ is merely uniformly continuous in $x$ was raised in \cite{DJK23,DK12}.

Recently, Liu \cite{Liu26} answered this question in the negative. He constructed a coefficient $a(x,y)$, which is symmetric and smooth in $y$ and whose modulus of continuity in $x$ is of order $(\log(1/r))^{-1/2}$, for which the interior $L_p$-estimates fail. This leaves a gap between the H\"older continuity assumed in \cite{DJK23} and the logarithmic modulus in \cite{Liu26}.

The purpose of this paper is to close this gap. We show that the right condition is the square-Dini condition
\begin{equation}\label{eq:sqD-intro}
\int_0^1\frac{\omega(t)^2}{t}\,dt<\infty
\end{equation}
on the modulus of continuity $\omega$ of $a(x,y)$ in $x$, uniformly in $y$. Under \eqref{eq:sqD-intro}, we prove that $\cL^a$ is bounded from $H^\sigma_{p,w}(\R^d)$ to $L_{p,w}(\R^d)$ for all $p\in(1,\infty)$ and $w\in A_p$ (Theorem~\ref{thm:bdd}), and that the global estimates and the unique solvability in weighted $L_p$ spaces, with mixed norms and weights in both $t$ and $x$ in the parabolic case, as well as the interior $L_p$-estimates, hold for the corresponding elliptic and parabolic equations (Theorems~\ref{thm:global}, \ref{thm:parabolic}, and \ref{thm:interior}). Conversely, by adapting the construction in \cite{Liu26}, we show that for every modulus $\omega$ violating \eqref{eq:sqD-intro} and satisfying a mild regularity condition, there is a coefficient with modulus of continuity $\omega$ for which the boundedness, the global estimates, and the interior estimates all fail (Theorem~\ref{thm:sharp}). In particular, for $\omega(t)=(\log(1/t))^{-\theta}$, the $L_p$-theory holds if and only if $\theta>1/2$ (Corollary~\ref{cor:log}). We note that the modulus $(\log(1/t))^{-1}$ satisfies \eqref{eq:sqD-intro}, but it is not a Dini modulus. This is in contrast with pointwise estimates, such as Schauder or $C^1$-type estimates, for which Dini-type conditions are natural; see, for instance, \cite{DK13,DK17}. Conditions of square-Dini type also appear in other contexts, for instance, in the study of elliptic measure; see \cite{FKP91}.

Our approach is different from that in \cite{DJK23}. There, the coefficient is represented by means of Riesz potentials of its fractional derivatives in $x$, which leads to the boundedness for $p>d/\beta$, where $\beta$ is the H\"older exponent, and then the Rubio de Francia extrapolation theorem is used to treat general $p\in(1,\infty)$. Instead, we use the Littlewood--Paley decomposition $u=\sum_{j\ge0}\Delta_ju$, and for each $j$ we decompose the coefficient in $x$ at frequency $2^{j-3}$. The key observation is that for a single Littlewood--Paley piece,
\[
\int_{\R^d}\big|\Delta_ju(x+y)-\Delta_ju(x)-(1_{\sigma>1}+1_{\sigma=1}1_{B_{2^{-j}}}(y))(\nabla\Delta_ju(x),y)\big|\frac{dy}{|y|^{d+\sigma}}\le N2^{j\sigma}\cM(\Delta_ju)(x),
\]
where $\cM$ is the Hardy--Littlewood maximal function. Hence any bounded measurable coefficient acts on each piece with a pointwise bound. The low-frequency part of the coefficient keeps the output frequency-localized near $2^j$, so that these outputs are almost orthogonal and only $\sup|a|$ enters. The high-frequency part of the coefficient has size $\omega(2^{-j})$, and we sum its contributions by the Cauchy--Schwarz inequality, which gives the factor $(\sum_j\omega(2^{-j})^2)^{1/2}$. This is the only place where the modulus of continuity enters. The argument works for each $p\in(1,\infty)$ directly and does not use extrapolation.

After completing this work, we learned of the recent paper \cite{Ryu26}, where the boundedness of non-local operators and weighted $L_p$-estimates for the corresponding elliptic equations were established under the stronger Dini continuity assumption in $x$.

The remaining part of the paper is organized as follows. In the next section, we introduce some notation and state the main results of the paper. In Section~\ref{sec:aux}, we present some auxiliary results. In Section~\ref{sec:bdd}, we prove the boundedness results. The $L_p$-estimates are proved in Section~\ref{sec:est}. Finally, in Section~\ref{sec:sharp}, we prove the sharpness of the square-Dini condition.

\medskip
\noindent\textbf{AI usage statement.} During the preparation of this work, the authors used Claude Opus~5.5. The authors first tried to find a counterexample as in \cite{Liu26} when the modulus of continuity of the coefficients in $x$ is of order $(\log(1/r))^{-1}$, which is almost Dini. This attempt indicated that the relevant threshold is the weaker square-Dini condition, under which a positive result should hold. The AI tool then proposed the Littlewood--Paley argument in the proof of Proposition~\ref{prop:main} and assisted in drafting the manuscript. The extension to weighted and mixed-norm spaces, the weaker ellipticity condition, and the treatment of the case $\sigma=1$ in the commutator estimates were provided by the authors, who also suggested several simplifications of the proofs. All arguments were checked and revised by the authors.

\section{Notation and main results}\label{sec:main}

For $x\in\R^d$ and $R\in(0,\infty)$, we denote
\[
B_R(x)=\{y\in\R^d:|x-y|<R\},\qquad B_R=B_R(0),
\]
where $d$ is a positive integer. As usual, $\cS(\R^d)$ is the Schwartz function space in $\R^d$, and $C_0^\infty(\R^d)$ is the set of all infinitely differentiable functions with compact support in $\R^d$. For $f\in\cS(\R^d)$, we use the Fourier transform $\hat f(\xi)=\int_{\R^d}e^{-i(x,\xi)}f(x)\,dx$. By $\partial^\sigma$ we mean the fractional Laplacian of order $\sigma$. That is,
\[
\partial^\sigma u(x)=-(-\Delta)^{\sigma/2}u(x)=\frac{c_{d,\sigma}}2\int_{\R^d}\big(u(x+y)+u(x-y)-2u(x)\big)\frac{dy}{|y|^{d+\sigma}}
\]
for sufficiently regular $u$ defined on $\R^d$, where
\[
c_{d,\sigma}=\Big(\int_{\R^d}(1-\cos y_1)\frac{dy}{|y|^{d+\sigma}}\Big)^{-1}.
\]

For $p\in(1,\infty)$ and a positive integer $k$, we denote by $A_p(\R^k)$ the set of all nonnegative functions $w$ on $\R^k$ such that
\[
[w]_{A_p(\R^k)}:=\sup_{x_0\in\R^k,R>0}\Big(\frac1{|B_R|}\int_{B_R(x_0)}w\,dx\Big)\Big(\frac1{|B_R|}\int_{B_R(x_0)}w^{-1/(p-1)}\,dx\Big)^{p-1}<\infty,
\]
where $B_R(x_0)$ is the ball in $\R^k$. When $k=d$, we simply write $[w]_{A_p}$.
For $w\in A_p(\R^d)$, by $L_{p,w}(\R^d)$ we mean the space of all measurable functions on $\R^d$ with the norm $\|f\|_{L_{p,w}}=(\int_{\R^d}|f|^pw\,dx)^{1/p}$. The weighted Bessel potential space is
\[
H^\sigma_{p,w}(\R^d)=\{u\in L_{p,w}(\R^d):(1-\Delta)^{\sigma/2}u\in L_{p,w}(\R^d)\},\qquad\|u\|_{H^\sigma_{p,w}}=\|(1-\Delta)^{\sigma/2}u\|_{L_{p,w}}.
\]
By \cite[Lemma~3.4]{DJK23},
\begin{equation}\label{eq:Hequiv}
\|u\|_{H^\sigma_{p,w}}\approx\|u\|_{L_{p,w}}+\|\partial^\sigma u\|_{L_{p,w}}.
\end{equation}
When $w=1$, we simply remove $w$ in the above notation. For $T\in(0,\infty)$, we write $L_p(T)=L_p((0,T)\times\R^d)$, and we denote by $\cH^\sigma_p(T)$ the space of $H^\sigma_p(\R^d)$-valued functions $u$ on $(0,T)$ with the norm $(\int_0^T\|u(t,\cdot)\|_{H^\sigma_p}^p\,dt)^{1/p}$. As solution spaces for parabolic equations, $\cH^{1,\sigma}_p(T)$ denotes the collection of functions $u\in\cH^\sigma_p(T)$ such that $u_t\in L_p(T)$, with the norm $\|u\|_{\cH^\sigma_p(T)}+\|u_t\|_{L_p(T)}$. For mixed norms with weights, let $p,q\in(1,\infty)$ and $w(t,x)=w_1(t)w_2(x)$ with $w_1\in A_p(\R)$ and $w_2\in A_q(\R^d)$. We denote by $L_{p,q,w}(T)$ the space of measurable functions $f$ on $(0,T)\times\R^d$ with the norm
\[
\|f\|_{L_{p,q,w}(T)}=\Big(\int_0^T\Big(\int_{\R^d}|f|^qw_2\,dx\Big)^{p/q}w_1\,dt\Big)^{1/p}.
\]
The spaces $\cH^\sigma_{p,q,w}(T)$ and $\cH^{1,\sigma}_{p,q,w}(T)$ are defined in the same way as above, with $H^\sigma_{q,w_2}(\R^d)$, $L_{p,q,w}(T)$, and $w_1(t)\,dt$ in place of $H^\sigma_p(\R^d)$, $L_p(T)$, and $dt$. Finally, $\cH^{1,\sigma}_{p,q,w,0}(T)$ denotes the closure in $\cH^{1,\sigma}_{p,q,w}(T)$ of the set of functions $u\in C^\infty([0,T]\times\R^d)$ vanishing for large $|x|$ and satisfying $u(0,\cdot)=0$. When $p=q$, we write $L_{p,w}(T)=L_{p,p,w}(T)$, $\cH^{1,\sigma}_{p,w,0}(T)=\cH^{1,\sigma}_{p,p,w,0}(T)$, and so on. We also use the weight
\[
\mu(x)=(1+|x|)^{-d-\sigma},\qquad\|u\|_{L_p(\mu)}=\Big(\int_{\R^d}|u|^p\mu\,dx\Big)^{1/p}.
\]

Throughout the paper, we denote
\[
\nabla^\sigma_yu(x)=u(x+y)-u(x)-\chi_\sigma(y)(\nabla u(x),y)
\]
for $x,y\in\R^d$, where $\chi_\sigma$ is as in \eqref{eq:L}. For a measurable function $c=c(x,y)$ on $\R^d\times\R^d$, which we call a coefficient, $\cL^c$ is defined as in \eqref{eq:L} with $c$ in place of $a$. If $c(x,-y)=c(x,y)$, then for every $\sigma\in(0,2)$,
\[
\cL^cu(x)=\frac12\int_{\R^d}\big(u(x+y)+u(x-y)-2u(x)\big)\,c(x,y)\,\frac{dy}{|y|^{d+\sigma}}.
\]
In particular, $\partial^\sigma=\cL^{c_{d,\sigma}}$. When $\sigma=1$, we impose the cancellation condition
\begin{equation}\label{eq:C}
\int_{r<|y|<R}y\,c(x,y)\,\frac{dy}{|y|^{d+1}}=0\qquad\text{for all }x\in\R^d\text{ and }0<r<R,
\end{equation}
which holds automatically if $c$ is symmetric in $y$.

Throughout the paper, we may omit $\R^d$ in $C_0^\infty(\R^d)$, $\cS(\R^d)$, or $L_p(\R^d)$ whenever the omission is clear from the context. Sometimes we use `sup' to represent the essential supremum. We write $N(d,\nu,\ldots)$ in the estimates to express that the constant $N$ is determined only by the parameters $d,\nu,\ldots$. The constant $N$ can differ from line to line.

We now state our assumption on the coefficient $a$. In the parabolic case, $a=a(t,x,y)$ is a measurable function on $(0,T)\times\R^d\times\R^d$, and in the elliptic case $a=a(x,y)$ is independent of $t$.

\begin{assumption}\label{ass}
\begin{enumerate}[label=(\roman*)]
\item There exist positive constants $\nu$ and $K_1$ and a nonnegative measurable function $m_0=m_0(t,y)$, which is homogeneous of degree zero in $y$, such that
\[
\inf_{\xi\in S^{d-1}}\int_{S^{d-1}}|(\xi,y)|^\sigma m_0(t,y)\,dS(y)\ge\nu\quad\text{for all }t
\]
and
\[
(2-\sigma)m_0(t,y)\le a(t,x,y)\le(2-\sigma)K_1\quad\text{for all }t,x,y.
\]
In the elliptic case, $m_0=m_0(y)$ is independent of $t$.
\item There is a nondecreasing function $\omega:[0,\infty)\to[0,\infty)$ such that
\[
|a(t,x_1,y)-a(t,x_2,y)|\le(2-\sigma)\,\omega(|x_1-x_2|)
\]
for all $t$ and $x_1,x_2,y\in\R^d$, and
\begin{equation}\label{eq:sqD}
[\omega]:=\Big(\int_0^1\frac{\omega(t)^2}{t}\,dt\Big)^{1/2}<\infty.
\end{equation}
\item If $\sigma=1$, then $a(t,\cdot,\cdot)$ satisfies \eqref{eq:C} for all $t$.
\end{enumerate}
\end{assumption}

\begin{remark}\label{rem:ell}
Assumption~\ref{ass} (i) is the ellipticity condition in \cite[Assumption~2.1]{DL23}; see also \cite{MP14,Zha13}. It holds, for instance, if $a\ge(2-\sigma)\delta$ for a constant $\delta>0$, in which case we may take $m_0\equiv\delta$. It also allows $a(t,x,\cdot)$ to vanish in a cone of directions. The lower bound is used only through the estimate for coefficients independent of $x$ in Lemma~\ref{lem:DK}.
\end{remark}

\begin{remark}\label{rem:sqD}
Since $\omega$ is nondecreasing, \eqref{eq:sqD} implies $\omega(0+)=0$. If $\omega$ is a Dini modulus, i.e., $\int_0^1\omega(t)t^{-1}\,dt<\infty$, then \eqref{eq:sqD} holds because $\omega(t)^2\le\omega(1)\omega(t)$ for $t\le1$. In particular, \eqref{eq:sqD} holds under the H\"older continuity assumption \cite[Assumption~2.1 (ii)]{DJK23}. The converse is not true: $\omega(t)=(\log(1/t))^{-\theta}$ near $t=0$ satisfies \eqref{eq:sqD} if and only if $\theta>1/2$, and it is a Dini modulus if and only if $\theta>1$.
\end{remark}

Our first main result is about the boundedness of non-local operators. Note that the lower bound in Assumption~\ref{ass} (i) is not needed here.

\begin{theorem}\label{thm:bdd}
Let $\sigma\in(0,2)$, $p\in(1,\infty)$, and $w\in A_p(\R^d)$. Let $a=a(x,y)$ be a measurable function satisfying $|a|\le(2-\sigma)K_1$ and Assumption~\ref{ass} (ii) and (iii). Then $\cL^a$ is a bounded operator from $H^\sigma_{p,w}(\R^d)$ to $L_{p,w}(\R^d)$. More precisely, there exists a constant $N=N(d,p,\sigma,[w]_{A_p})$ such that, for any $u\in H^\sigma_{p,w}(\R^d)$,
\begin{equation}\label{eq:bdd}
\|\cL^au\|_{L_{p,w}}\le N(2-\sigma)\Big((K_1+[\omega])\|\partial^\sigma u\|_{L_{p,w}}+K_1\|u\|_{L_{p,w}}\Big).
\end{equation}
\end{theorem}

\begin{remark}\label{rem:const}
We do not pursue the dependence of the constants on $\sigma$, for instance, as $\sigma\nearrow2$; compare \cite{DJK23}. On the other hand, the proof of Theorem~\ref{thm:bdd} uses no regularity of $a$ in $y$. The cancellation condition \eqref{eq:C} is used only when $\sigma=1$.
\end{remark}

When $a=a(t,x,y)$, by using Theorem~\ref{thm:bdd} for each $t\in(0,T)$ and then integrating with respect to $t$, we obtain the following corollary for the parabolic operator.

\begin{corollary}\label{cor:par-bdd}
Let $\sigma\in(0,2)$, $p,q\in(1,\infty)$, $T\in(0,\infty)$, and $w(t,x)=w_1(t)w_2(x)$ with $w_1\in A_p(\R)$ and $w_2\in A_q(\R^d)$. Suppose that $a=a(t,x,y)$ satisfies $|a|\le(2-\sigma)K_1$ and Assumption~\ref{ass} (ii) and (iii). Then the operator $\partial_t-\cL^a$ is bounded from $\cH^{1,\sigma}_{p,q,w}(T)$ to $L_{p,q,w}(T)$.
\end{corollary}

As the second set of our main results, we prove the following weighted $L_p$-estimates and the unique solvability for elliptic and parabolic equations.

\begin{theorem}\label{thm:global}
Let $\sigma\in(0,2)$, $p\in(1,\infty)$, $w\in A_p(\R^d)$, and let $a=a(x,y)$ satisfy Assumption~\ref{ass}. Then there exist constants $\lambda_0\ge1$, depending only on $d,p,\sigma,\nu,K_1,\omega,[w]_{A_p}$, and $N$, depending only on $d,p,\sigma,\nu,K_1,[w]_{A_p}$, such that, for any $\lambda\ge\lambda_0$ and $u\in H^\sigma_{p,w}(\R^d)$,
\begin{equation}\label{eq:global}
\|\partial^\sigma u\|_{L_{p,w}}+\lambda\|u\|_{L_{p,w}}\le N\|\lambda u-\cL^au\|_{L_{p,w}}.
\end{equation}
Moreover, for any $\lambda\ge\lambda_0$ and $f\in L_{p,w}(\R^d)$, there exists a unique $u\in H^\sigma_{p,w}(\R^d)$ satisfying \eqref{eq:ell}.
\end{theorem}

\begin{theorem}\label{thm:parabolic}
Let $\sigma\in(0,2)$, $p,q\in(1,\infty)$, $T\in(0,\infty)$, $w(t,x)=w_1(t)w_2(x)$ with $w_1\in A_p(\R)$ and $w_2\in A_q(\R^d)$, and let $a=a(t,x,y)$ satisfy Assumption~\ref{ass}. Then there exist constants $\lambda_0\ge1$, depending only on $d$, $p$, $q$, $\sigma$, $\nu$, $K_1$, $\omega$, $[w_1]_{A_p(\R)}$, and $[w_2]_{A_q}$, and $N$, depending only on $d$, $p$, $q$, $\sigma$, $\nu$, $K_1$, $[w_1]_{A_p(\R)}$, and $[w_2]_{A_q}$, such that the following holds. For any $\lambda\ge\lambda_0$, if $u\in\cH^{1,\sigma}_{p,q,w,0}(T)$ satisfies
\begin{equation}\label{eq:par}
u_t-\cL^au+\lambda u=f\quad\text{in }(0,T)\times\R^d
\end{equation}
with $f\in L_{p,q,w}(T)$, then
\begin{equation}\label{eq:par-est}
\|u_t\|_{L_{p,q,w}(T)}+\|\partial^\sigma u\|_{L_{p,q,w}(T)}+\lambda\|u\|_{L_{p,q,w}(T)}\le N\|f\|_{L_{p,q,w}(T)}.
\end{equation}
Moreover, for any $\lambda\ge\lambda_0$ and $f\in L_{p,q,w}(T)$, there exists a unique $u\in\cH^{1,\sigma}_{p,q,w,0}(T)$ satisfying \eqref{eq:par}.
\end{theorem}

We also have the following interior estimates, which are of the type considered in \cite[Section~7]{DK12} and \cite{Liu26}.

\begin{theorem}\label{thm:interior}
Let $\sigma\in(0,2)$ and $p\in(1,\infty)$. Let $a:B_4\times\R^d\to\R$ be measurable and satisfy Assumption~\ref{ass} for $x,x_1,x_2\in B_4$. Then, for any $u\in C^2(\R^d)$ with $\|u\|_{L_p(\mu)}<\infty$, the following hold.
\begin{enumerate}[label=(\roman*)]
\item If $\sigma\in(0,1)$, then
\[
\|\partial^\sigma u\|_{L_p(B_1)}\le N\big(\|\cL^au\|_{L_p(B_2)}+\|u\|_{L_p(\mu)}\big).
\]
\item If $\sigma\in[1,2)$, then for any $\theta\in(0,1]$,
\[
\|\partial^\sigma u\|_{L_p(B_1)}\le N\|\cL^au\|_{L_p(B_2)}+\theta\|Du\|_{L_p(B_4)}+N_\theta\|u\|_{L_p(\mu)}.
\]
\end{enumerate}
Here $N=N(d,p,\sigma,\nu,K_1,\omega)$ and $N_\theta=N_\theta(d,p,\sigma,\nu,K_1,\omega,\theta)$.
\end{theorem}

Our last main result shows that the square-Dini condition \eqref{eq:sqD} is sharp. We consider moduli $\omega:[0,\infty)\to[0,\infty)$ satisfying
\begin{equation}\label{eq:notsqD}
\begin{gathered}
\omega\text{ is continuous and nondecreasing},\quad\omega(0)=0,\quad\omega>0\text{ on }(0,\infty),\\
t\mapsto\omega(t)/t\text{ is nonincreasing on }(0,1/2],\quad\text{and}\quad\int_0^{1/2}\frac{\omega(t)^2}{t}\,dt=\infty.
\end{gathered}
\end{equation}

\begin{theorem}\label{thm:sharp}
Let $\sigma\in(0,2)$, and let $\omega$ satisfy \eqref{eq:notsqD}. Then there exist a coefficient $a\in C^\infty(\R^d\times(\R^d\setminus\{0\}))$, which is symmetric in $y$, and a constant $C\ge1$ such that
\begin{equation}\label{eq:sharp-a}
C^{-1}\le a\le C,\qquad|a(x_1,y)-a(x_2,y)|\le C\,\omega(|x_1-x_2|),\qquad\big|D^i_y\big(a(x,y)|y|^{-d-\sigma}\big)\big|\le C|y|^{-d-\sigma-i}
\end{equation}
for all $x,x_1,x_2\in\R^d$, $y\ne0$, and $i=1,2$, and such that for every $p\in(1,\infty)$ the following hold.
\begin{enumerate}[label=(\alph*)]
\item The operator $\cL^a$ is not bounded from $H^\sigma_p(\R^d)$ to $L_p(\R^d)$.
\item For any $\lambda\ge0$ and $N>0$, the estimate \eqref{eq:global} with $w=1$ fails for some $u\in C_0^\infty(\R^d)$.
\item The interior estimates in Theorem~\ref{thm:interior} fail. More precisely, if $\sigma\ne1$, there is no constant $N$ such that
\[
\|\partial^\sigma u\|_{L_p(B_1)}\le N\big(\|\cL^au\|_{L_p(B_2)}+1_{\sigma>1}\|Du\|_{L_p(B_4)}+\|u\|_{L_p(\mu)}\big)
\]
for all $u\in C^\infty(\R^d)$ whose derivatives of all orders are bounded. If $\sigma=1$, there is $\theta\in(0,1]$ such that no constant $N_\theta$ satisfies
\[
\|\partial^1u\|_{L_p(B_1)}\le N_\theta\|\cL^au\|_{L_p(B_2)}+\theta\|Du\|_{L_p(B_4)}+N_\theta\|u\|_{L_p(\mu)}
\]
for all such $u$.
\end{enumerate}
\end{theorem}

By \eqref{eq:sharp-a}, the coefficient in Theorem~\ref{thm:sharp} satisfies Assumption~\ref{ass} (i) and (iii), and (ii) with $\omega$ replaced by $C\omega/(2-\sigma)$, except for \eqref{eq:sqD}. It is also smooth in $y$ and satisfies the derivative bounds in $y$ assumed in \cite{Liu26}.

\begin{remark}\label{rem:regular}
The assumption that $\omega(t)/t$ is nonincreasing in \eqref{eq:notsqD} holds for any concave $\omega$ with $\omega(0)=0$. Moreover, if $\omega$ is subadditive, as is the case for the optimal modulus of continuity of a uniformly continuous function, then it is well known that its least concave majorant $\bar\omega$ satisfies $\omega\le\bar\omega\le2\omega$. Hence, if a subadditive $\omega$ satisfies $\int_0^{1/2}\omega(t)^2t^{-1}\,dt=\infty$, then Theorem~\ref{thm:sharp} applied to $\bar\omega$ gives a coefficient whose modulus of continuity is bounded by $2C\omega$.
\end{remark}

\begin{corollary}\label{cor:log}
For $\theta>0$, let $\omega_\theta(0)=0$, $\omega_\theta(t)=(\log(1/t))^{-\theta}$ for $0<t\le1/2$, and $\omega_\theta(t)=(\log2)^{-\theta}$ for $t>1/2$.
\begin{enumerate}[label=(\roman*)]
\item If $\theta>1/2$, then Theorems~\ref{thm:bdd}--\ref{thm:interior} hold with $\omega=\omega_\theta$. This includes $\theta=1$, for which $\omega_1$ is not a Dini modulus.
\item If $0<\theta\le1/2$, then the conclusions of Theorem~\ref{thm:sharp} hold with $\omega=\omega_\theta$.
\end{enumerate}
\end{corollary}

\begin{proof}
We have $[\omega_\theta]^2\le N+\int_0^{1/2}t^{-1}(\log(1/t))^{-2\theta}\,dt$, which is finite if and only if $2\theta>1$. For (ii), note that on $(0,1/2]$,
\[
\frac{d}{dt}\log\frac{\omega_\theta(t)}{t}=\frac1t\Big(\frac{\theta}{\log(1/t)}-1\Big)<0,
\]
because $\log(1/t)\ge\log2>\theta$. Hence $\omega_\theta$ satisfies \eqref{eq:notsqD}, and Theorem~\ref{thm:sharp} applies.
\end{proof}

\begin{remark}\label{rem:liu}
Liu \cite{Liu26} considers the operator $\frac12\int(u(x+y)+u(x-y)-2u(x))K(x,y)\,dy$ with a symmetric kernel $K$. With $a(x,y)=|y|^{d+\sigma}K(x,y)$, this operator is $\cL^a$, and the assumptions in \cite[(1.2), (1.4)]{Liu26} imply Assumption~\ref{ass} (i) and (ii) on $B_4$, up to the normalizing factor $2-\sigma$ and without \eqref{eq:sqD}. The derivative bounds in $y$ in \cite[(1.3)]{Liu26} are not used in Theorems~\ref{thm:bdd}--\ref{thm:interior}. Note also that the interior estimates in Theorem~\ref{thm:interior} (ii) hold with an arbitrarily small constant in front of $\|Du\|_{L_p(B_4)}$, not only for $\sigma=1$ but also for $\sigma\in(1,2)$; cf. \cite[(7.2), (7.3)]{DK12}.
\end{remark}

\begin{remark}\label{rem:weights}
Theorem~\ref{thm:interior} is stated in unweighted spaces. In its proof, the non-local tail of $u$ is estimated by the Minkowski inequality and measured in $L_p(\mu)$, which does not extend directly to $L_{p,w}$ with a general $A_p$ weight. We do not pursue weighted interior estimates here.
\end{remark}

\section{Auxiliary results}\label{sec:aux}

In this section, we present some auxiliary results. Except for Lemma~\ref{lem:DK}, all functions and coefficients in this section are independent of the time variable.

\subsection{Maximal functions}
For $f\in L_{1,\text{loc}}(\R^d)$, let $\cM f(x)=\sup_{r>0}|B_r|^{-1}\int_{B_r(x)}|f|\,dy$ be the Hardy--Littlewood maximal function of $f$. We shall use the following well-known facts.

\begin{lemma}\label{lem:max}
\begin{enumerate}[label=(\roman*)]
\item Let $\Psi\in L_1(\R^d)$ satisfy $|\Psi(x)|\le\Psi^*(|x|)$, where $\Psi^*$ is nonincreasing and $\int_{\R^d}\Psi^*(|x|)\,dx<\infty$. Then $|\Psi_t*f|\le\cM f\int_{\R^d}\Psi^*(|x|)\,dx$ for all $t>0$, where $\Psi_t=t^{-d}\Psi(\cdot/t)$.
\item For $\gamma\in(0,2)$, $R>0$, and $f\in L_{1,\text{loc}}(\R^d)$,
\[
\int_{|y|\ge R}|f(x+y)|\frac{dy}{|y|^{d+\gamma}}\le N\gamma^{-1}R^{-\gamma}\cM f(x),\qquad\int_{|y|<R}|f(x+y)|\frac{dy}{|y|^{d-\gamma}}\le N\gamma^{-1}R^{\gamma}\cM f(x),
\]
where $N=N(d)$.
\item Let $p\in(1,\infty)$ and $w\in A_p(\R^d)$. Then for any sequence $\{f_j\}$,
\[
\Big\|\Big(\sum_j(\cM f_j)^2\Big)^{1/2}\Big\|_{L_{p,w}}\le N\Big\|\Big(\sum_j|f_j|^2\Big)^{1/2}\Big\|_{L_{p,w}},
\]
where $N=N(d,p,[w]_{A_p})$.
\end{enumerate}
\end{lemma}

\begin{proof}
For (i), see \cite{Ste70}. Assertion (ii) is \cite[Lemma~3.2]{DJK23}. Assertion (iii) is the Fefferman--Stein inequality \cite{FS71} and its weighted version \cite{AJ81}. In particular, (iii) contains the weighted Hardy--Littlewood maximal function theorem.
\end{proof}

\subsection{Littlewood--Paley decomposition}
Fix a radial $\Phi\in\cS(\R^d)$ such that $\hat\Phi=1$ on $B_1$ and $\supp\hat\Phi\subset B_2$, and let $\hat\psi(\xi)=\hat\Phi(\xi)-\hat\Phi(2\xi)$. Then $\supp\hat\psi\subset\{1/2\le|\xi|\le2\}$. For $u\in\cS(\R^d)$, we set
\[
\Delta_0u=\Phi*u,\qquad\Delta_ju=\cF^{-1}\big[\hat\psi(2^{-j}\xi)\hat u\big],\quad j\ge1,
\]
and we write $u_j=\Delta_ju$ for $j\ge0$.
For $j\ge1$, $\hat u_j$ is supported in $\{2^{j-1}\le|\xi|\le2^{j+1}\}$. Since
\[
\sum_{j=0}^Ju_j=\cF^{-1}\big[\hat\Phi(2^{-J}\xi)\hat u\big],
\]
we have
\begin{equation}\label{eq:LP}
u=\sum_{j\ge0}u_j\quad\text{in }\cS(\R^d).
\end{equation}
To decompose coefficients in $x$, let $\Theta\in\cS(\R^d)$ be given by $\hat\Theta(\xi)=\hat\Phi(2\xi)$. Then $\int_{\R^d}\Theta\,dx=1$ and $\supp\hat\Theta\subset B_1$. For $k\in\mathbb{Z}$, let $\Theta_k=2^{kd}\Theta(2^k\cdot)$, and for a bounded measurable coefficient $c(x,y)$, we denote by
\begin{equation}\label{eq:Tk}
\Theta_k*c(x,y)=\int_{\R^d}\Theta_k(z)\,c(x-z,y)\,dz
\end{equation}
the convolution of $c$ and $\Theta_k$ in $x$. For each fixed $y$, the distributional Fourier transform of $x\mapsto\Theta_k*c(x,y)$ is supported in $\overline{B_{2^k}}$.

\begin{lemma}\label{lem:LP}
Let $u\in\cS(\R^d)$.
\begin{enumerate}[label=(\roman*)]
\item For $j\ge0$, $k=0,1,2$, and $x\in\R^d$,
\[
\sup_{|z-x|\le2^{-j}}|D^ku_j(z)|\le N2^{kj}\cM(u_j)(x),
\]
where $N=N(d)$.
\item Let $p\in(1,\infty)$ and $w\in A_p(\R^d)$. Then
\[
\Big\|\Big(\sum_{j\ge1}4^{j\sigma}|u_j|^2\Big)^{1/2}\Big\|_{L_{p,w}}\le N\|\partial^\sigma u\|_{L_{p,w}},
\]
where $N=N(d,p,\sigma,[w]_{A_p})$.
\item Let $p\in(1,\infty)$ and $w\in A_p(\R^d)$. If $F_j\in L_\infty(\R^d)$, $j\ge1$, and $\hat F_j$ is supported in $\{2^{j-2}\le|\xi|\le2^{j+2}\}$, then for any integers $1\le J_1\le J_2$,
\[
\Big\|\sum_{j=J_1}^{J_2}F_j\Big\|_{L_{p,w}}\le N\Big\|\Big(\sum_{j=J_1}^{J_2}|F_j|^2\Big)^{1/2}\Big\|_{L_{p,w}},
\]
where $N=N(d,p,[w]_{A_p})$.
\end{enumerate}
\end{lemma}

\begin{proof}
(i) The function $\hat u_j$ is supported in $B_{2^{j+1}}$. Fix $\Upsilon\in\cS(\R^d)$ with $\hat\Upsilon=1$ on $B_2$, and set $\Upsilon_j=2^{jd}\Upsilon(2^j\cdot)$. Then $u_j=\Upsilon_j*u_j$, so that
\[
D^ku_j(z)=\int_{\R^d}D^k\Upsilon_j(z-y)\,u_j(y)\,dy,\qquad|D^k\Upsilon_j(z-y)|\le N2^{j(d+k)}(1+2^j|z-y|)^{-d-1}.
\]
If $|z-x|\le2^{-j}$, then $1+2^j|x-y|\le2(1+2^j|z-y|)$. Hence
\[
|D^ku_j(z)|\le N2^{kj}\int_{\R^d}2^{jd}(1+2^j|x-y|)^{-d-1}|u_j(y)|\,dy\le N2^{kj}\cM(u_j)(x)
\]
by Lemma~\ref{lem:max} (i).

(ii) We have $2^{j\sigma}u_j=\Psi^{(j)}*\partial^\sigma u$, where $\hat\Psi(\xi)=-\hat\psi(\xi)|\xi|^{-\sigma}$ belongs to $C_0^\infty(\R^d\setminus\{0\})$ and $\Psi^{(j)}=2^{jd}\Psi(2^j\cdot)$. Thus the assertion follows from the weighted Littlewood--Paley inequality; see, for instance, \cite{Kur80}.

(iii) Let $\Upsilon'\in\cS(\R^d)$ satisfy $\hat\Upsilon'=1$ on $\{1/4\le|\xi|\le4\}$ and $\supp\hat\Upsilon'\subset\{1/8\le|\xi|\le8\}$, and let $\Upsilon'_j=2^{jd}\Upsilon'(2^j\cdot)$. Then $F_j=\Upsilon'_j*F_j$. For $g\in C_0^\infty(\R^d)$, by H\"older's inequality,
\[
\Big|\int_{\R^d}\sum_{j=J_1}^{J_2}F_j\,g\,dx\Big|\le\Big\|\Big(\sum_{j=J_1}^{J_2}|F_j|^2\Big)^{1/2}\Big\|_{L_{p,w}}\Big\|\Big(\sum_{j\ge1}|\tilde\Upsilon'_j*g|^2\Big)^{1/2}\Big\|_{L_{p',w'}},
\]
where $\tilde\Upsilon'_j(x)=\Upsilon'_j(-x)$, $p'=p/(p-1)$, and $w'=w^{-1/(p-1)}\in A_{p'}$. By the weighted Littlewood--Paley inequality, the last factor is bounded by $N\|g\|_{L_{p',w'}}$. The assertion follows by duality.
\end{proof}

\subsection{Moduli of continuity}
For a bounded measurable coefficient $c$, let
\[
\vp_c(t)=\sup\big\{|c(x+h,y)-c(x,y)|:\ x,y\in\R^d,\ |h|\le t\big\}.
\]
This function is nondecreasing and subadditive. Hence
\begin{equation}\label{eq:subadd}
\vp_c(\lambda t)\le(1+\lambda)\vp_c(t),\qquad\lambda,t>0.
\end{equation}
If $\vp$ is a nondecreasing function, then
\[
\vp(2^{-j})^2\log2\le\int_{2^{-j}}^{2^{1-j}}\frac{\vp(t)^2}{t}\,dt,\qquad j\ge1,
\]
and thus
\begin{equation}\label{eq:sum-int}
\sum_{j\ge1}\vp(2^{-j})^2\le(\log2)^{-1}\int_0^1\frac{\vp(t)^2}{t}\,dt.
\end{equation}

\subsection{Further auxiliary lemmas}

The following lemma is about operators with coefficients independent of $x$.

\begin{lemma}\label{lem:DK}
Let $\sigma\in(0,2)$, $p\in(1,\infty)$, $T\in(0,\infty)$, and $w(t,x)=w_1(t)w_2(x)$ with $w_1\in A_p(\R)$ and $w_2\in A_p(\R^d)$. Let $a_0=a_0(t,y)$ be measurable and satisfy $(2-\sigma)m_0(t,y)\le a_0(t,y)\le(2-\sigma)K_1$, where $m_0$ and $\nu$ are as in Assumption~\ref{ass} (i), and \eqref{eq:C} for each $t$ if $\sigma=1$. Then for any $\lambda\ge0$ and $u\in\cH^{1,\sigma}_{p,w,0}(T)$,
\[
\|u_t\|_{L_{p,w}(T)}+\|\partial^\sigma u\|_{L_{p,w}(T)}+\lambda\|u\|_{L_{p,w}(T)}\le N_0\|u_t-\cL^{a_0}u+\lambda u\|_{L_{p,w}(T)},
\]
where $N_0=N_0(d,p,\sigma,\nu,K_1,[w_1]_{A_p(\R)},[w_2]_{A_p})$ is independent of $T$.
\end{lemma}

\begin{proof}
This is contained in \cite[Theorem~2.9]{DL23} with $p=q$. With the kernel $a_0(t,y)|y|^{-d-\sigma}$, the conditions on $a_0$ are those in \cite[Assumption~2.1]{DL23} with $\Lambda=K_1$. When $w=1$, see also \cite[Proposition~1]{MP14}. In \cite{DL23}, the cancellation condition is imposed on spheres when $\sigma=1$. It follows from \eqref{eq:C} for almost every radius, which suffices for the proof there.
\end{proof}

Finally, we consider commutators with cut-off functions. For $f\in C_0^\infty(\R^d)$ and $g\in C^2(\R^d)$ with $\|g\|_{L_1(\mu)}<\infty$, a direct computation shows that
\begin{equation}\label{eq:comm}
\cL^c(fg)=f\,\cL^cg+g\,\cL^cf+I^c[f,g],
\end{equation}
where
\[
I^c[f,g](x)=\int_{\R^d}\big(f(x+y)-f(x)\big)\big(g(x+y)-g(x)\big)\,c(x,y)\,\frac{dy}{|y|^{d+\sigma}}.
\]
Indeed, the terms involving $\chi_\sigma$ cancel identically, so \eqref{eq:comm} holds for every $\sigma\in(0,2)$, and all the integrals converge absolutely. Since $\partial^\sigma=\cL^{c_{d,\sigma}}$, \eqref{eq:comm} holds with $\partial^\sigma$ in place of $\cL^c$ and $I^{c_{d,\sigma}}$ in place of $I^c$.

\begin{lemma}\label{lem:avcomm}
Let $\sigma\in(0,2)$, $p\in(1,\infty)$, $w\in A_p(\R^d)$, and let $c$ be a measurable function with $|c|\le(2-\sigma)K_1$. Let $\zeta\in C_0^\infty(\R^d)$ with $\|\zeta\|_{L_p}=1$, and set $\zeta_z=\zeta(\cdot-z)$ for $z\in\R^d$. Then for any $u\in C_0^\infty(\R^d)$,
\begin{equation}\label{eq:avcomm1}
\Big(\int_{\R^d}\|u\,\cL^c\zeta_z\|_{L_{p,w}}^p\,dz\Big)^{1/p}\le N\|u\|_{L_{p,w}}.
\end{equation}
Moreover, if $\sigma\in(0,1)$, then
\[
\Big(\int_{\R^d}\|I^c[\zeta_z,u]\|_{L_{p,w}}^p\,dz\Big)^{1/p}\le N\|u\|_{L_{p,w}},
\]
and if $\sigma\in[1,2)$, then for any $\theta\in(0,1]$,
\[
\Big(\int_{\R^d}\|I^c[\zeta_z,u]\|_{L_{p,w}}^p\,dz\Big)^{1/p}\le\theta\|Du\|_{L_{p,w}}+N_\theta\|u\|_{L_{p,w}}.
\]
Here $N=N(d,p,\sigma,K_1,\zeta,[w]_{A_p})$, and $N_\theta$ depends in addition on $\theta$.
\end{lemma}

\begin{proof}
By Fubini's theorem, for any functions $F_z(x)$,
\begin{equation}\label{eq:fubini-z}
\int_{\R^d}\|F_z\|_{L_{p,w}}^p\,dz=\int_{\R^d}\Big(\int_{\R^d}|F_z(x)|^p\,dz\Big)w(x)\,dx.
\end{equation}

Let $g(x)=\int_{\R^d}|\nabla^\sigma_y\zeta(x)|\,|y|^{-d-\sigma}\,dy$. Then $g$ is bounded. If $\supp\zeta\subset B_R$ and $|x|\ge2R$, then $\zeta(x)=0$ and $\nabla\zeta(x)=0$, and $\zeta(x+y)\neq0$ implies $|y|\ge|x|/2$. Thus $g(x)\le\|\zeta\|_{L_1}(|x|/2)^{-d-\sigma}$, so that $g\in L_p$. Since $|\cL^c\zeta_z(x)|\le(2-\sigma)K_1g(x-z)$, \eqref{eq:avcomm1} follows from \eqref{eq:fubini-z}.

For the estimates of $I^c[\zeta_z,u]$, we fix $x$ and use the Minkowski inequality in $L_p(dz)$:
\[
\Big(\int_{\R^d}|I^c[\zeta_z,u](x)|^p\,dz\Big)^{1/p}\le(2-\sigma)K_1\int_{\R^d}\|\zeta(\cdot+y)-\zeta\|_{L_p}|u(x+y)-u(x)|\frac{dy}{|y|^{d+\sigma}}\le NJ(x),
\]
where we used $\|\zeta(\cdot+y)-\zeta\|_{L_p}\le\min\{|y|\|D\zeta\|_{L_p},2\}$ and
\[
J(x)=\int_{\R^d}\min\{|y|,1\}|u(x+y)-u(x)|\frac{dy}{|y|^{d+\sigma}}.
\]
By \eqref{eq:fubini-z}, it suffices to estimate $\|J\|_{L_{p,w}}$. Since $u$ is continuous, $|u|\le\cM u$. By Lemma~\ref{lem:max} (ii),
\[
\int_{|y|\ge1}|u(x+y)-u(x)|\frac{dy}{|y|^{d+\sigma}}\le N\cM u(x).
\]
If $\sigma<1$, again by Lemma~\ref{lem:max} (ii),
\[
\int_{|y|<1}|u(x+y)-u(x)|\frac{dy}{|y|^{d+\sigma-1}}\le N\cM u(x).
\]
If $\sigma\ge1$, we fix $r\in(0,1)$. For $|y|<r$, we use $|u(x+y)-u(x)|\le|y|\int_0^1|Du(x+sy)|\,ds$ and Lemma~\ref{lem:max} (ii) to get
\[
\begin{aligned}
\int_{|y|<r}|u(x+y)-u(x)|\frac{dy}{|y|^{d+\sigma-1}}&\le\int_0^1s^{\sigma-2}\int_{|z|<sr}|Du(x+z)|\frac{dz}{|z|^{d+\sigma-2}}\,ds\\
&\le N(2-\sigma)^{-1}r^{2-\sigma}\cM(Du)(x).
\end{aligned}
\]
For $r\le|y|<1$, we have $|y|^{1-d-\sigma}\le r^{-\sigma}|y|^{1-d}$, and by Lemma~\ref{lem:max} (ii),
\[
\int_{r\le|y|<1}|u(x+y)-u(x)|\frac{dy}{|y|^{d+\sigma-1}}\le r^{-\sigma}\int_{|y|<1}\big(|u(x+y)|+|u(x)|\big)\frac{dy}{|y|^{d-1}}\le Nr^{-\sigma}\cM u(x).
\]
Combining these estimates with the weighted Hardy--Littlewood maximal function theorem, we obtain $\|J\|_{L_{p,w}}\le N\|u\|_{L_{p,w}}$ if $\sigma<1$, and
\[
\|J\|_{L_{p,w}}\le N(2-\sigma)^{-1}r^{2-\sigma}\|Du\|_{L_{p,w}}+Nr^{-\sigma}\|u\|_{L_{p,w}}
\]
if $\sigma\ge1$. It remains to choose $r$ sufficiently small.
\end{proof}

\begin{remark}\label{rem:avcomm}
The proof of Lemma~\ref{lem:avcomm} follows that of \cite[Lemma~3.6]{DJK23} and the estimate of $I_2$ in the proof of \cite[Theorem~2.3]{DJK23}, except that for $\sigma\ge1$ we split the integral at $|y|=r$ with a small $r$ instead of $|y|=1$. This gives a small constant in front of $\|Du\|_{L_{p,w}}$, which is needed when $\sigma=1$ for the absorption in Section~\ref{sec:est}. A similar argument is used in the proof of \cite[Corollary~2.10]{DL23}.
\end{remark}

\begin{lemma}\label{lem:cutoff}
Let $\sigma\in(0,2)$, $p\in(1,\infty)$, $R\ge1$, and let $c$ be a measurable function with $|c|\le(2-\sigma)K_1$. Let $\eta\in C_0^\infty(B_R)$ with $0\le\eta\le1$, and let $u\in C^2(\R^d)$ with $\|u\|_{L_p(\mu)}<\infty$.
\begin{enumerate}[label=(\roman*)]
\item If $\sigma\in(0,1)$, then
\[
\|u\,\cL^c\eta\|_{L_p}+\|I^c[\eta,u]\|_{L_p}\le N\|u\|_{L_p(\mu)}.
\]
\item If $\sigma\in[1,2)$, then for any $\theta\in(0,1]$,
\[
\|u\,\cL^c\eta\|_{L_p}+\|I^c[\eta,u]\|_{L_p}\le\theta\|Du\|_{L_p(B_{R+2})}+N_\theta\|u\|_{L_p(\mu)}.
\]
\end{enumerate}
Here $N=N(d,p,\sigma,K_1,R,\eta)$ and $N_\theta=N_\theta(d,p,\sigma,K_1,R,\eta,\theta)$. No regularity of $c$ in $x$ and no cancellation condition are needed.
\end{lemma}

\begin{proof}
Since $\mu\ge(R+3)^{-d-\sigma}$ on $B_{R+2}$, we have $\|u\|_{L_p(B_{R+2})}\le N\|u\|_{L_p(\mu)}$. Also, $\|u\|_{L_1(\mu)}\le N\|u\|_{L_p(\mu)}$ by H\"older's inequality.

\emph{Estimate of $u\,\cL^c\eta$.} For $|x|\le2R$, we have $|\cL^c\eta(x)|\le N$. For $|x|>2R$, $\eta$ and $\nabla\eta$ vanish at $x$, and $\eta(x+y)\neq0$ implies $|y|\ge|x|/2$, so that $|\cL^c\eta(x)|\le N|x|^{-d-\sigma}$. Hence $|u\,\cL^c\eta|\le N|u|\mu$, and $\|u\,\cL^c\eta\|_{L_p}^p\le N\int|u|^p\mu^p\,dx\le N\|u\|_{L_p(\mu)}^p$.

\emph{Estimate of $I^c[\eta,u]$ for $|y|\ge1$.} Here we use
\[
\big|\big(\eta(x+y)-\eta(x)\big)\big(u(x+y)-u(x)\big)\big|\le\big(\eta(x+y)+\eta(x)\big)\big(|u(x+y)|+|u(x)|\big),
\]
which gives four terms. First, by the Minkowski inequality, the $L_p$ norm of $\int_{|y|\ge1}\eta(x+y)|u(x+y)||y|^{-d-\sigma}dy$ is bounded by $N\|u\|_{L_p(B_R)}$. Second, by the same reasoning as for $\cL^c\eta$,
\[
|u(x)|\int_{|y|\ge1}\eta(x+y)\frac{dy}{|y|^{d+\sigma}}\le N|u(x)|\mu(x).
\]
Third, for $x\in B_R$ and $|y|\ge1$, we have $1+|x+y|\le(R+2)|y|$, so that $|y|^{-d-\sigma}\le(R+2)^{d+\sigma}\mu(x+y)$. Thus
\[
\eta(x)\int_{|y|\ge1}|u(x+y)|\frac{dy}{|y|^{d+\sigma}}\le N1_{B_R}(x)\|u\|_{L_1(\mu)}.
\]
Fourth, $\eta(x)|u(x)|\int_{|y|\ge1}|y|^{-d-\sigma}dy\le N1_{B_R}(x)|u(x)|$. All four terms have $L_p$ norms bounded by $N\|u\|_{L_p(\mu)}$.

\emph{Estimate of $I^c[\eta,u]$ for $|y|<1$.} Here $|\eta(x+y)-\eta(x)|\le\|D\eta\|_{L_\infty}|y|1_{B_{R+1}}(x)$. If $\sigma<1$, this part is bounded by
\[
N1_{B_{R+1}}(x)\int_{|y|<1}\big(|u(x+y)|+|u(x)|\big)\frac{dy}{|y|^{d+\sigma-1}},
\]
whose $L_p$ norm is bounded by $N\|u\|_{L_p(B_{R+2})}$ by the Minkowski inequality. If $\sigma\ge1$, we fix $\rho\in(0,1)$. For $|y|<\rho$, we use $|u(x+y)-u(x)|\le|y|\int_0^1|Du(x+sy)|\,ds$, and note that $x+sy\in B_{R+2}$ if $x\in B_{R+1}$. By the Minkowski inequality in $(s,y)$, the $L_p$ norm of this part is bounded by
\[
N\int_{|y|<\rho}\frac{dy}{|y|^{d+\sigma-2}}\,\|Du\|_{L_p(B_{R+2})}=N(2-\sigma)^{-1}\rho^{2-\sigma}\|Du\|_{L_p(B_{R+2})}.
\]
For $\rho\le|y|<1$, we use $|u(x+y)-u(x)|\le|u(x+y)|+|u(x)|$, and the $L_p$ norm of this part is bounded by $N(\rho)\|u\|_{L_p(B_{R+2})}$. It remains to choose $\rho$ so that $N(2-\sigma)^{-1}\rho^{2-\sigma}\le\theta$.
\end{proof}

\section{Boundedness}\label{sec:bdd}

In this section, we prove Theorem~\ref{thm:bdd} and Corollary~\ref{cor:par-bdd}. The key step is Proposition~\ref{prop:main} below, which holds for any bounded measurable coefficient and gives an explicit dependence on its modulus of continuity in $x$.

For $\sigma=1$ and $j\ge0$, we set
\[
\nabla^{1,j}_yv(x)=v(x+y)-v(x)-1_{B_{2^{-j}}}(y)(\nabla v(x),y),
\]
and for $\sigma\neq1$, we set $\nabla^{\sigma,j}_y=\nabla^\sigma_y$. Recall that $u_j=\Delta_ju$. By Lemma~\ref{lem:LP} (i), for $u\in\cS(\R^d)$,
\begin{equation}\label{eq:locsup}
\sup_{|z-x|\le2^{-j}}|D^ku_j(z)|\le N2^{kj}\cM(u_j)(x),\qquad k=0,1,2,\ j\ge0.
\end{equation}
Moreover, since $|u_0|\le N\cM u$ by Lemma~\ref{lem:max} (i), for $p\in(1,\infty)$ and $w\in A_p(\R^d)$, we have
\begin{equation}\label{eq:m0}
\|\cM(u_0)\|_{L_{p,w}}\le N\|u_0\|_{L_{p,w}}\le N\|u\|_{L_{p,w}}.
\end{equation}

\begin{lemma}\label{lem:block}
For $u\in\cS(\R^d)$, $j\ge0$, and $x\in\R^d$, we have
\begin{equation}\label{eq:block}
\int_{\R^d}|\nabla^{\sigma,j}_yu_j(x)|\frac{dy}{|y|^{d+\sigma}}\le N2^{j\sigma}\cM(u_j)(x),
\end{equation}
where $N=N(d,\sigma)$. Consequently, if $c=c(x,y)$ is a measurable function satisfying $|c|\le K$, and \eqref{eq:C} when $\sigma=1$, then
\begin{equation}\label{eq:block2}
|\cL^cu_j(x)|\le NK\,2^{j\sigma}\cM(u_j)(x).
\end{equation}
\end{lemma}

\begin{proof}
We split the integral at $|y|=2^{-j}$.

\emph{The case $|y|<2^{-j}$.} If $\sigma<1$, then by \eqref{eq:locsup}, $|u_j(x+y)-u_j(x)|\le|y|\sup_{B_{2^{-j}}(x)}|Du_j|\le N2^j|y|\cM(u_j)(x)$. Hence this part is bounded by
\[
N2^j\cM(u_j)(x)\int_{|y|<2^{-j}}\frac{dy}{|y|^{d+\sigma-1}}=N(1-\sigma)^{-1}2^{j\sigma}\cM(u_j)(x).
\]
If $\sigma\ge1$, then by \eqref{eq:locsup}, for $|y|<2^{-j}$,
\[
|\nabla^{\sigma,j}_yu_j(x)|\le|y|^2\sup_{B_{2^{-j}}(x)}|D^2u_j|\le N2^{2j}|y|^2\cM(u_j)(x).
\]
Hence this part is bounded by $N(2-\sigma)^{-1}2^{j\sigma}\cM(u_j)(x)$.

\emph{The case $|y|\ge2^{-j}$.} Here
\[
|\nabla^{\sigma,j}_yu_j(x)|\le|u_j(x+y)|+|u_j(x)|+1_{\sigma>1}|y||\nabla u_j(x)|.
\]
By Lemma~\ref{lem:max} (ii) with $\gamma=\sigma$,
\[
\int_{|y|\ge2^{-j}}\big(|u_j(x+y)|+|u_j(x)|\big)\frac{dy}{|y|^{d+\sigma}}\le N\sigma^{-1}2^{j\sigma}\big(\cM(u_j)(x)+|u_j(x)|\big)\le N2^{j\sigma}\cM(u_j)(x),
\]
because $|u_j|\le\cM(u_j)$ for the continuous function $u_j$. If $\sigma>1$, by \eqref{eq:locsup},
\[
|\nabla u_j(x)|\int_{|y|\ge2^{-j}}\frac{dy}{|y|^{d+\sigma-1}}=N(\sigma-1)^{-1}2^{j(\sigma-1)}|\nabla u_j(x)|\le N2^{j\sigma}\cM(u_j)(x).
\]
This proves \eqref{eq:block}.

To prove \eqref{eq:block2} when $\sigma\neq1$, we use $|\cL^cu_j(x)|\le K\int|\nabla^\sigma_yu_j(x)||y|^{-d-\sigma}dy$. If $\sigma=1$ and $j\ge1$, then by \eqref{eq:C} for $c$,
\[
\int_{2^{-j}\le|y|<1}(\nabla u_j(x),y)\,c(x,y)\frac{dy}{|y|^{d+1}}=0,
\]
so that $\cL^cu_j(x)=\int\nabla^{1,j}_yu_j(x)\,c(x,y)|y|^{-d-1}dy$, and \eqref{eq:block} applies. For $j=0$, there is nothing to change. The lemma is proved.
\end{proof}

\begin{proposition}\label{prop:main}
Let $\sigma\in(0,2)$, $p\in(1,\infty)$, and $w\in A_p(\R^d)$. Let $c=c(x,y)$ be a measurable function satisfying $|c|\le K$, and \eqref{eq:C} if $\sigma=1$. Then for any $u\in\cS(\R^d)$,
\begin{equation}\label{eq:main}
\|\cL^cu\|_{L_{p,w}}\le N\Big(K+\Big(\sum_{j\ge1}\min\{\vp_c(2^{-j}),K\}^2\Big)^{1/2}\Big)\|\partial^\sigma u\|_{L_{p,w}}+NK\|u\|_{L_{p,w}},
\end{equation}
where $N=N(d,p,\sigma,[w]_{A_p})$.
\end{proposition}

\begin{proof}
\emph{Step 1: decomposition of the coefficient.} For $j\ge1$, we write $c=\Theta_{j-3}*c+r_j$, where $\Theta_{j-3}*c$ is defined in \eqref{eq:Tk}. Then $|\Theta_{j-3}*c|\le\|\Theta\|_{L_1}K$. Moreover, since $\int\Theta=1$, by \eqref{eq:subadd} and the change of variables $z=2^{3-j}z'$,
\begin{equation}\label{eq:rj}
\begin{aligned}
|r_j(x,y)|&=\Big|\int_{\R^d}\Theta_{j-3}(z)\big(c(x,y)-c(x-z,y)\big)\,dz\Big|\\
&\le\int_{\R^d}|\Theta(z')|\min\big\{(1+8|z'|)\vp_c(2^{-j}),2K\big\}\,dz'\le N\min\{\vp_c(2^{-j}),K\}.
\end{aligned}
\end{equation}
Since the convolution is only in $x$, both $\Theta_{j-3}*c$ and $r_j$ satisfy \eqref{eq:C} when $\sigma=1$.

\emph{Step 2: decomposition of the operator.} For $v\in C^2(\R^d)$ with bounded derivatives up to order two, we have $|\cL^cv(x)|\le NK(\|v\|_{L_\infty}+\|Dv\|_{L_\infty}+\|D^2v\|_{L_\infty})$. By \eqref{eq:LP}, the partial sums of $\sum_{j\ge0}u_j$ converge to $u$ in this norm. Hence for every $x\in\R^d$,
\begin{equation}\label{eq:split}
\cL^cu(x)=\cL^cu_0(x)+\lim_{J\to\infty}\sum_{j=1}^J\big(G_j(x)+H_j(x)\big),
\end{equation}
where $G_j=\cL^{\Theta_{j-3}*c}u_j$ and $H_j=\cL^{r_j}u_j$.

\emph{Step 3: estimate of $\cL^cu_0$.} By \eqref{eq:block2} with $j=0$, we have $|\cL^cu_0|\le NK\cM(u_0)$. Thus, by \eqref{eq:m0}, $\|\cL^cu_0\|_{L_{p,w}}\le NK\|u\|_{L_{p,w}}$.

\emph{Step 4: estimate of $\sum_jG_j$.} By \eqref{eq:block2} applied to the coefficient $\Theta_{j-3}*c$, $|G_j|\le NK2^{j\sigma}\cM(u_j)$, so that $G_j\in L_\infty$. We claim that $\hat G_j$ is supported in $\{2^{j-2}\le|\xi|\le2^{j+2}\}$. Fix $j\ge 1$ and $y$, and let
\[
g_y(x)=\nabla^{\sigma,j}_yu_j(x)\cdot \Theta_{j-3}*c(x,y),
\]
where, when $\sigma=1$, the truncation $1_{B_{2^{-j}}}$ is used as in the proof of Lemma~\ref{lem:block}. The first factor is a Schwartz function whose Fourier transform is supported in $\{2^{j-1}\le|\xi|\le2^{j+1}\}$, and the second factor is a bounded smooth function whose Fourier transform is supported in $\overline{B_{2^{j-3}}}$. Thus $\hat g_y$ is supported in $\{2^{j-2}\le|\xi|\le2^{j+2}\}$. Now let $h\in\cS(\R^d)$ be such that $\hat h$ vanishes in a neighborhood of this set. By \eqref{eq:block},
\[
\int_{\R^d}\int_{\R^d}|g_y(x)||h(x)|\frac{dy}{|y|^{d+\sigma}}\,dx\le NK\int_{\R^d}2^{j\sigma}\cM(u_j)|h|\,dx<\infty.
\]
Hence, by Fubini's theorem, $\int G_jh\,dx=\int(\int g_yh\,dx)|y|^{-d-\sigma}dy=0$, which proves the claim. By Lemma~\ref{lem:LP} (iii), Lemma~\ref{lem:max} (iii), and Lemma~\ref{lem:LP} (ii), for any $1\le J_1\le J_2$,
\[
\Big\|\sum_{j=J_1}^{J_2}G_j\Big\|_{L_{p,w}}\le NK\Big\|\Big(\sum_{j=J_1}^{J_2}4^{j\sigma}(\cM(u_j))^2\Big)^{1/2}\Big\|_{L_{p,w}},
\]
and
\begin{equation}\label{eq:sqfn}
\Big\|\Big(\sum_{j\ge1}4^{j\sigma}(\cM(u_j))^2\Big)^{1/2}\Big\|_{L_{p,w}}\le N\|\partial^\sigma u\|_{L_{p,w}}.
\end{equation}
By the dominated convergence theorem, the partial sums of $\sum_jG_j$ form a Cauchy sequence in $L_{p,w}$, and their limit $G$ satisfies $\|G\|_{L_{p,w}}\le NK\|\partial^\sigma u\|_{L_{p,w}}$.

\emph{Step 5: estimate of $\sum_jH_j$.} By \eqref{eq:block2} applied to the coefficient $r_j$, together with \eqref{eq:rj}, $|H_j|\le N\min\{\vp_c(2^{-j}),K\}2^{j\sigma}\cM(u_j)$. By the Cauchy--Schwarz inequality,
\[
\sum_{j\ge1}|H_j|\le N\Big(\sum_{j\ge1}\min\{\vp_c(2^{-j}),K\}^2\Big)^{1/2}\Big(\sum_{j\ge1}4^{j\sigma}(\cM(u_j))^2\Big)^{1/2},
\]
and, by \eqref{eq:sqfn},
\[
\Big\|\sum_{j\ge1}|H_j|\Big\|_{L_{p,w}}\le N\Big(\sum_{j\ge1}\min\{\vp_c(2^{-j}),K\}^2\Big)^{1/2}\|\partial^\sigma u\|_{L_{p,w}}.
\]
We note that this is the only step where the modulus of continuity of $c$ is used.

\emph{Step 6: conclusion.} By \eqref{eq:split} and Step 5, the partial sums $\sum_{j\le J}G_j(x)$ converge for almost every $x$, and their limit coincides almost everywhere with the function $G$ in Step 4. Therefore, $|\cL^cu|\le|\cL^cu_0|+|G|+\sum_{j\ge1}|H_j|$ almost everywhere, and \eqref{eq:main} follows from Steps 3--5.
\end{proof}

\begin{proof}[Proof of Theorem~\ref{thm:bdd}]
We apply Proposition~\ref{prop:main} with $c=a$ and $K=(2-\sigma)K_1$. By Assumption~\ref{ass} (ii) and \eqref{eq:sum-int}, $\sum_{j\ge1}\vp_a(2^{-j})^2\le(2-\sigma)^2(\log2)^{-1}[\omega]^2$. This gives \eqref{eq:bdd} for $u\in\cS(\R^d)$. Since $\cS(\R^d)$ is dense in $H^\sigma_{p,w}(\R^d)$, the theorem follows from \eqref{eq:Hequiv}.
\end{proof}

\begin{proof}[Proof of Corollary~\ref{cor:par-bdd}]
For $u\in\cH^{1,\sigma}_{p,q,w}(T)$, by Theorem~\ref{thm:bdd} applied to $a(t,\cdot,\cdot)$, the exponent $q$, and the weight $w_2$ for each $t\in(0,T)$, we have $\|\cL^au(t,\cdot)\|_{L_{q,w_2}}\le N\|u(t,\cdot)\|_{H^\sigma_{q,w_2}}$ with $N$ independent of $t$. Taking the $L_p$ norm in $t$ with respect to $w_1(t)\,dt$ gives $\|\cL^au\|_{L_{p,q,w}(T)}\le N\|u\|_{\cH^\sigma_{p,q,w}(T)}$, which proves the corollary. Note that the argument works for any nonnegative function $w_1$.
\end{proof}

\section{\texorpdfstring{$L_p$}{Lp}-estimates}\label{sec:est}

In this section, we prove Theorems~\ref{thm:global}, \ref{thm:parabolic}, and \ref{thm:interior} by the freezing coefficient argument; see, for instance, the proof of \cite[Lemma~1.6.3]{Kry08} for second-order equations and the proof of \cite[Lemma~5.2]{DJK23} for non-local equations. The following lemma replaces \cite[Lemma~5.1]{DJK23}.

\begin{lemma}\label{lem:loc}
Let $\sigma\in(0,2)$, $p\in(1,\infty)$, $w\in A_p(\R^d)$, and let $a=a(x,y)$ be a measurable function satisfying $|a|\le(2-\sigma)K_1$ and Assumption~\ref{ass} (ii) and (iii). Let $\varepsilon\in(0,1/8]$, and let $\eta\in C_0^\infty(B_{4\varepsilon})$ satisfy $0\le\eta\le1$, $\eta=1$ on $B_{2\varepsilon}$, and $|D\eta|\le1/\varepsilon$. Set $a_0(y)=a(0,y)$ and
\[
\varrho(\varepsilon)=\omega(4\varepsilon)+\Big(\sum_{j\ge1:\,2^{-j}<\varepsilon}\omega(2^{-j})^2\Big)^{1/2}+\omega(8\varepsilon)\big(1+\log_2(1/\varepsilon)\big)^{1/2}.
\]
Then for any $u\in\cS(\R^d)$,
\begin{equation}\label{eq:loc}
\|\eta(\cL^a-\cL^{a_0})u\|_{L_{p,w}}\le N(2-\sigma)\Big(\varrho(\varepsilon)\|\partial^\sigma u\|_{L_{p,w}}+\omega(4\varepsilon)\|u\|_{L_{p,w}}\Big),
\end{equation}
where $N=N(d,p,\sigma,[w]_{A_p})$. Moreover, $\varrho(\varepsilon)\to0$ as $\varepsilon\to0$.
\end{lemma}

\begin{proof}
Let $c(x,y)=\eta(x)\big(a(x,y)-a(0,y)\big)$. Since $\eta$ multiplies the coefficient at the point $x$, we have $\eta(\cL^a-\cL^{a_0})u=\cL^cu$. Moreover, $c$ satisfies \eqref{eq:C} when $\sigma=1$, and $|c|\le(2-\sigma)\omega(4\varepsilon)=:K$.

We estimate $\vp_c$. For any $t>0$, since $|c|\le K$,
\[
\vp_c(t)\le2K=2(2-\sigma)\omega(4\varepsilon).
\]
Now let $t<\varepsilon$ and $|x-x'|\le t$. If neither $x$ nor $x'$ belongs to $B_{4\varepsilon}$, then $c(x,y)=c(x',y)=0$. Otherwise, both points belong to $B_{5\varepsilon}\subset B_{8\varepsilon}$, so that $|a(x',y)-a(0,y)|\le(2-\sigma)\omega(8\varepsilon)$ and
\begin{align*}
|c(x,y)-c(x',y)|&\le\eta(x)|a(x,y)-a(x',y)|+|a(x',y)-a(0,y)||\eta(x)-\eta(x')|\\
&\le(2-\sigma)\big(\omega(t)+\omega(8\varepsilon)t/\varepsilon\big).
\end{align*}
Hence $\vp_c(t)\le(2-\sigma)(\omega(t)+\omega(8\varepsilon)t/\varepsilon)$ for $t<\varepsilon$. There are at most $\log_2(1/\varepsilon)$ integers $j\ge1$ with $2^{-j}\ge\varepsilon$, and $\sum_{2^{-j}<\varepsilon}(2^{-j}/\varepsilon)^2\le4/3$. Using the first bound for $2^{-j}\ge\varepsilon$ and the second one for $2^{-j}<\varepsilon$, we obtain
\[
\sum_{j\ge1}\vp_c(2^{-j})^2\le(2-\sigma)^2\Big(2\sum_{2^{-j}<\varepsilon}\omega(2^{-j})^2+\frac83\omega(8\varepsilon)^2+4\omega(4\varepsilon)^2\log_2\frac1\varepsilon\Big)\le4(2-\sigma)^2\varrho(\varepsilon)^2.
\]
Since also $K\le(2-\sigma)\varrho(\varepsilon)$, Proposition~\ref{prop:main} gives \eqref{eq:loc}.

It remains to show that $\varrho(\varepsilon)\to0$. By \eqref{eq:sqD} and \eqref{eq:sum-int}, the second term in $\varrho(\varepsilon)$ is the tail of a convergent series, and $\omega(0+)=0$. Since $\omega$ is nondecreasing,
\[
\frac12\omega(8\varepsilon)^2\log\frac1{8\varepsilon}\le\int_{8\varepsilon}^{\sqrt{8\varepsilon}}\frac{\omega(t)^2}{t}\,dt\to0\quad\text{as }\varepsilon\to0,
\]
so that $\omega(8\varepsilon)^2\log(1/\varepsilon)\to0$. The lemma is proved.
\end{proof}

\begin{lemma}\label{lem:small}
Let $\sigma\in(0,2)$, $p\in(1,\infty)$, and $w(t,x)=w_1(t)w_2(x)$ with $w_1\in A_p(\R)$ and $w_2\in A_p(\R^d)$. There exist $\varepsilon_0\in(0,1/8]$ and $N'$, depending only on $d,p,\sigma,\nu,K_1,\omega,[w_1]_{A_p(\R)},[w_2]_{A_p}$, such that the following holds. Let $T\in(0,\infty)$, and let $a=a(t,x,y)$ satisfy Assumption~\ref{ass}. Then for any $x_0\in\R^d$, $\lambda\ge0$, and $u\in\cH^{1,\sigma}_{p,w,0}(T)$ such that $u(t,\cdot)$ vanishes outside $B_{\varepsilon_0}(x_0)$ for every $t\in(0,T)$,
\[
\|u_t\|_{L_{p,w}(T)}+\|\partial^\sigma u\|_{L_{p,w}(T)}+\lambda\|u\|_{L_{p,w}(T)}\le2N_0\|u_t-\cL^au+\lambda u\|_{L_{p,w}(T)}+N'\|u\|_{L_{p,w}(T)},
\]
where $N_0$ is the constant in Lemma~\ref{lem:DK}.
\end{lemma}

\begin{proof}
By a standard approximation argument, we may assume that $u$ is smooth. By translation in $x$, we may assume $x_0=0$; note that a translation of an $A_p$ weight is an $A_p$ weight with the same constant. Let $\varepsilon\in(0,1/8]$, assume that $u(t,\cdot)\in C_0^\infty(B_\varepsilon)$ for each $t$, and let $\eta$ be as in Lemma~\ref{lem:loc}. Set $a_0(t,y)=a(t,0,y)$. Since $a_0$ satisfies the conditions of Lemma~\ref{lem:DK}, we have
\begin{multline*}
\|u_t\|_{L_{p,w}(T)}+\|\partial^\sigma u\|_{L_{p,w}(T)}+\lambda\|u\|_{L_{p,w}(T)}\le N_0\Big(\|u_t-\cL^au+\lambda u\|_{L_{p,w}(T)}\\
+\|\eta(\cL^a-\cL^{a_0})u\|_{L_{p,w}(T)}+\|(1-\eta)(\cL^a-\cL^{a_0})u\|_{L_{p,w}(T)}\Big).
\end{multline*}
For each $t\in(0,T)$, the second term is estimated by Lemma~\ref{lem:loc} applied to $a(t,\cdot,\cdot)$, $u(t,\cdot)$, and the weight $w_2$. For the third term, let $|x|\ge2\varepsilon$. Then $u(t,x)=0$ and $\nabla u(t,x)=0$, so that
\[
(\cL^a-\cL^{a_0})u(t,x)=\int_{|y|\ge\varepsilon}u(t,x+y)\big(a(t,x,y)-a(t,0,y)\big)\frac{dy}{|y|^{d+\sigma}},
\]
because $u(t,x+y)\neq0$ implies $|y|\ge\varepsilon$. By Lemma~\ref{lem:max} (ii), the absolute value of the right-hand side is bounded by $NK_1\varepsilon^{-\sigma}\cM u(t,\cdot)(x)$. By the weighted Hardy--Littlewood maximal function theorem with the weight $w_2$, and integrating in $t$ with respect to $w_1(t)\,dt$, we obtain
\begin{multline*}
\|\eta(\cL^a-\cL^{a_0})u\|_{L_{p,w}(T)}+\|(1-\eta)(\cL^a-\cL^{a_0})u\|_{L_{p,w}(T)}\\
\le N(2-\sigma)\varrho(\varepsilon)\|\partial^\sigma u\|_{L_{p,w}(T)}+N\big(\omega(4\varepsilon)+K_1\varepsilon^{-\sigma}\big)\|u\|_{L_{p,w}(T)},
\end{multline*}
where $N=N(d,p,\sigma,[w_2]_{A_p})$. We now fix $\varepsilon=\varepsilon_0$ so small that $N_0N(2-\sigma)\varrho(\varepsilon_0)\le1/2$, which is possible because $\varrho(\varepsilon)\to0$. Absorbing the term containing $\|\partial^\sigma u\|_{L_{p,w}(T)}$ into the left-hand side, we obtain the lemma.
\end{proof}

\begin{proof}[Proof of Theorem~\ref{thm:parabolic}]
\emph{Step 1: the a priori estimate when $p=q$.} We use a version of the partition of unity argument in the proof of \cite[Theorem~1.6.4]{Kry08}. By Corollary~\ref{cor:par-bdd} and the definition of $\cH^{1,\sigma}_{p,w,0}(T)$, we may assume that $u$ is smooth, $u(0,\cdot)=0$, and $u(t,\cdot)$ vanishes outside a fixed ball for all $t$. Let $\varepsilon_0$ and $N'$ be the constants in Lemma~\ref{lem:small}. Fix $\zeta\in C_0^\infty(B_{\varepsilon_0})$ with $\zeta\ge0$ and $\|\zeta\|_{L_p}=1$, and set $\zeta_z=\zeta(\cdot-z)$. Then for any $g\in L_{p,w}(T)$,
\begin{equation}\label{eq:pou}
\int_{\R^d}\|\zeta_zg\|_{L_{p,w}(T)}^p\,dz=\|g\|_{L_{p,w}(T)}^p.
\end{equation}
Since $\zeta_zu(t,\cdot)$ vanishes outside $B_{\varepsilon_0}(z)$, by Lemma~\ref{lem:small},
\begin{multline}\label{eq:zeta-est}
\|(\zeta_zu)_t\|_{L_{p,w}(T)}+\|\partial^\sigma(\zeta_zu)\|_{L_{p,w}(T)}+\lambda\|\zeta_zu\|_{L_{p,w}(T)}\\
\le2N_0\|(\zeta_zu)_t-\cL^a(\zeta_zu)+\lambda\zeta_zu\|_{L_{p,w}(T)}+N'\|\zeta_zu\|_{L_{p,w}(T)}.
\end{multline}
Since $\zeta_z$ is independent of $t$, by \eqref{eq:comm} we have
\[
(\zeta_zu)_t-\cL^a(\zeta_zu)+\lambda\zeta_zu=\zeta_zf-u\,\cL^a\zeta_z-I^a[\zeta_z,u]
\]
and
\[
\zeta_z\partial^\sigma u=\partial^\sigma(\zeta_zu)-u\,\partial^\sigma\zeta_z-I^{c_{d,\sigma}}[\zeta_z,u].
\]
We take the $L_p(dz)$ norms of both sides of \eqref{eq:zeta-est} and of the last identity, apply Lemma~\ref{lem:avcomm} with the weight $w_2$ for each $t\in(0,T)$, and integrate in $t$ with respect to $w_1(t)\,dt$. By \eqref{eq:pou} and the Minkowski inequality, for any $\theta\in(0,1]$,
\[
\|u_t\|_{L_{p,w}(T)}+\|\partial^\sigma u\|_{L_{p,w}(T)}+\lambda\|u\|_{L_{p,w}(T)}\le6N_0\|f\|_{L_{p,w}(T)}+1_{\sigma\ge1}\theta\|Du\|_{L_{p,w}(T)}+N_\theta\|u\|_{L_{p,w}(T)},
\]
where $N_\theta$ depends only on $d,p,\sigma,\nu,K_1,\omega,[w_1]_{A_p(\R)},[w_2]_{A_p}$, and $\theta$. If $\sigma\ge1$, then by \cite[Lemmas~3.4 and 3.5]{DJK23}, for each $t\in(0,T)$,
\[
\|Du(t,\cdot)\|_{L_{p,w_2}}\le N_1\big(\|\partial^\sigma u(t,\cdot)\|_{L_{p,w_2}}+\|u(t,\cdot)\|_{L_{p,w_2}}\big),
\]
where $N_1=N_1(d,p,\sigma,[w_2]_{A_p})$. Taking $\theta=1/(2N_1)$, we obtain
\[
\|u_t\|_{L_{p,w}(T)}+\|\partial^\sigma u\|_{L_{p,w}(T)}+\lambda\|u\|_{L_{p,w}(T)}\le12N_0\|f\|_{L_{p,w}(T)}+N''\|u\|_{L_{p,w}(T)}.
\]
For $\lambda\ge\lambda_0:=\max\{1,2N''\}$, the last term is bounded by $\frac\lambda2\|u\|_{L_{p,w}(T)}$, and we obtain \eqref{eq:par-est} with $N=24N_0$ when $p=q$. Note that $\lambda_0$ and $N$ are independent of $T$.

\emph{Step 2: the a priori estimate for general $p,q$.} We use the extrapolation theorem; see \cite[Theorem~2.5]{DK18} and also \cite{CMP11}. By the definition of $\cH^{1,\sigma}_{p,q,w,0}(T)$ and Corollary~\ref{cor:par-bdd}, it suffices to prove \eqref{eq:par-est} for $u\in C^\infty([0,T]\times\R^d)$ vanishing for large $|x|$ and satisfying $u(0,\cdot)=0$. For such $u$, we extend
\[
|u_t|+|\partial^\sigma u|+\lambda|u|\qquad\text{and}\qquad|u_t-\cL^au+\lambda u|
\]
by zero for $t\notin(0,T)$. By Step 1 with, for instance, $p=q=2$, for any $K_0\ge1$, there exist $\lambda_0$, depending only on $K_0,d,\sigma,\nu,K_1,\omega$, and $N$, depending only on $K_0,d,\sigma,\nu,K_1$, such that for any $\lambda\ge\lambda_0$, the $L_{2,w}$ norm of the first function on $\R\times\R^d$ is bounded by $N$ times that of the second function, for all weights $w=w_1(t)w_2(x)$ with $[w_1]_{A_2(\R)},[w_2]_{A_2}\le K_0$. By the extrapolation theorem, the corresponding estimate holds in $L_{p,q,w}$ for any $p,q\in(1,\infty)$ and $w=w_1(t)w_2(x)$ with $w_1\in A_p(\R)$ and $w_2\in A_q(\R^d)$, where $K_0$, and thus $\lambda_0$ and $N$, are determined by $d,p,q,[w_1]_{A_p(\R)}$, and $[w_2]_{A_q}$. This gives \eqref{eq:par-est}.

\emph{Step 3: solvability.} The uniqueness follows from \eqref{eq:par-est}. For the existence, for $\tau\in[0,1]$ let $a_\tau=\tau a+(1-\tau)(2-\sigma)K_1$. Each $a_\tau$ satisfies Assumption~\ref{ass} with the same $K_1$ and $\omega$, with $m_0$ replaced by $\tau m_0+(1-\tau)K_1$, and with $\nu$ replaced by
\[
\min\Big\{\nu,K_1\int_{S^{d-1}}|y_1|^\sigma\,dS(y)\Big\},
\]
so that \eqref{eq:par-est} holds for $\cL^{a_\tau}$ uniformly in $\tau$. By Corollary~\ref{cor:par-bdd} and the method of continuity, it suffices to consider $\tau=0$, that is, the equation with $\cL^a$ replaced by a positive multiple of $\partial^\sigma$, whose unique solvability in $\cH^{1,\sigma}_{p,q,w,0}(T)$ is contained in \cite[Theorem~2.9]{DL23}.
\end{proof}

\begin{proof}[Proof of Theorem~\ref{thm:global}]
We derive the a priori estimate from Theorem~\ref{thm:parabolic}; compare the proof of \cite[Corollary~2.9]{DJK23} and \cite[Remark~2.7]{DL23}. Let $\lambda_0$ and $N$ be the constants in Theorem~\ref{thm:parabolic}, which are independent of $T$, and let $\lambda\ge\lambda_0$. By Theorem~\ref{thm:bdd} and the density of $C_0^\infty(\R^d)$ in $H^\sigma_{p,w}(\R^d)$, we may assume $u\in C_0^\infty(\R^d)$. Let $f=\lambda u-\cL^au$. For $n\ge1$, the function $v(t,x)=(t/n)u(x)$ is smooth, vanishes for large $|x|$ and at $t=0$, and satisfies
\[
v_t-\cL^av+\lambda v=\frac tnf+\frac1nu\quad\text{in }(0,n)\times\R^d.
\]
Since $\big(\int_0^n(t/n)^p\,dt\big)^{1/p}=(n/(p+1))^{1/p}$, Theorem~\ref{thm:parabolic} with $T=n$, $w_1\equiv1$, and $w_2=w$ gives
\[
\Big(\frac n{p+1}\Big)^{1/p}\big(\|\partial^\sigma u\|_{L_{p,w}}+\lambda\|u\|_{L_{p,w}}\big)\le N\Big(\frac n{p+1}\Big)^{1/p}\|f\|_{L_{p,w}}+Nn^{1/p-1}\|u\|_{L_{p,w}}.
\]
Dividing by $(n/(p+1))^{1/p}$ and letting $n\to\infty$, we obtain \eqref{eq:global}.

For the solvability, let $a_\tau=\tau a+(1-\tau)(2-\sigma)K_1$, $\tau\in[0,1]$. As in the proof of Theorem~\ref{thm:parabolic}, \eqref{eq:global} holds for $\cL^{a_\tau}$ uniformly in $\tau$, so that by Theorem~\ref{thm:bdd} and the method of continuity, the solvability follows from that for $\tau=0$, that is, for the equation $\lambda u+c(-\Delta)^{\sigma/2}u=f$ with a positive constant $c$, which follows from the weighted Fourier multiplier theorem; see \cite[Lemma~3.3]{DJK23} and \cite{KW79}.
\end{proof}

\begin{proof}[Proof of Theorem~\ref{thm:interior}]
Let $\pi(x)=x$ for $|x|\le3$ and $\pi(x)=3x/|x|$ for $|x|>3$. Since $\pi$ is $1$-Lipschitz, the function $\tilde a(x,y)=a(\pi(x),y)$ satisfies Assumption~\ref{ass} in $\R^d$ with the same $m_0,\nu,K_1,\omega$, and $\cL^{\tilde a}v=\cL^av$ in $B_3$.

Let $\eta\in C_0^\infty(B_2)$ satisfy $0\le\eta\le1$ and $\eta=1$ on $B_1$. Then $\eta u\in C_0^2(\R^d)\subset H^\sigma_p(\R^d)$, and for such functions the extension of $\cL^{\tilde a}$ given by Theorem~\ref{thm:bdd} coincides with the pointwise definition \eqref{eq:L}. By Theorem~\ref{thm:global} with $w=1$ and $\lambda=\lambda_0$, and by \eqref{eq:comm},
\[
\|\partial^\sigma(\eta u)\|_{L_p}\le N\Big(\lambda_0\|u\|_{L_p(B_2)}+\|\eta\,\cL^{\tilde a}u\|_{L_p}+\|u\,\cL^{\tilde a}\eta\|_{L_p}+\|I^{\tilde a}[\eta,u]\|_{L_p}\Big).
\]
Since $\eta\,\cL^{\tilde a}u=\eta\,\cL^au$, we have $\|\eta\,\cL^{\tilde a}u\|_{L_p}\le\|\cL^au\|_{L_p(B_2)}$. Also, $\|u\|_{L_p(B_2)}\le N\|u\|_{L_p(\mu)}$. On the other hand, by \eqref{eq:comm} with $\partial^\sigma$ in place of $\cL^c$,
\[
\|\partial^\sigma u\|_{L_p(B_1)}\le\|\eta\,\partial^\sigma u\|_{L_p}\le\|\partial^\sigma(\eta u)\|_{L_p}+\|u\,\partial^\sigma\eta\|_{L_p}+\|I^{c_{d,\sigma}}[\eta,u]\|_{L_p}.
\]
It remains to apply Lemma~\ref{lem:cutoff} with $R=2$ to $\tilde a$ and to the constant coefficient $c_{d,\sigma}$. In the case $\sigma\in[1,2)$, we apply it with $\theta/(N+1)$ in place of $\theta$.
\end{proof}

\section{Sharpness of the square-Dini condition}\label{sec:sharp}

In this section, we prove Theorem~\ref{thm:sharp}. The construction is due to Liu \cite{Liu26}, who treated the case $\omega(t)=(\log(1/t))^{-1/2}$. We present it for a general modulus satisfying \eqref{eq:notsqD}, and for the reader's convenience we give all the details.

Throughout this section, $\sigma\in(0,2)$ and $\omega$ satisfies \eqref{eq:notsqD}. All coefficients are symmetric in $y$, so that
\[
\cL^au(x)=\frac12\int_{\R^d}\delta_yu(x)\,a(x,y)\frac{dy}{|y|^{d+\sigma}},\qquad\delta_yu(x)=u(x+y)+u(x-y)-2u(x).
\]
For $k\ge0$, let $b_k(x)=\cos(2^kx_1)$. By the definition of $\partial^\sigma$ and the change of variables $y\mapsto2^{-k}y$,
\begin{equation}\label{eq:frac-ek}
(-\Delta)^{\sigma/2}b_k=c_{d,\sigma}\,b_k\int_{\R^d}\big(1-\cos(2^ky_1)\big)\frac{dy}{|y|^{d+\sigma}}=2^{k\sigma}b_k.
\end{equation}

\subsection{Building blocks}
Fix a nonnegative radial $\varphi\in C_0^\infty(\R^d)$ with $\supp\varphi\subset\{3/4<|\xi|<5/4\}$ and $\varphi(e_1)=1$, where $e_1=(1,0,\ldots,0)$. Let $q\in\cS(\R^d)$ be given by $\hat q(\xi)=|\xi|^\sigma\varphi(\xi)$, and for $j\ge1$ let
\[
q_j(y)=2^{j(d+\sigma)}q(2^jy),\qquad\text{so that}\qquad\hat q_j(\xi)=|\xi|^\sigma\varphi(2^{-j}\xi).
\]
Since $\hat q$ is real, radial, smooth, and compactly supported, $q$ is a real radial Schwartz function, and $\int q=\hat q(0)=0$. The same holds for each $q_j$.

\begin{lemma}\label{lem:qsum}
For $i=0,1,2$, there is a constant $C_i=C_i(d,\sigma,\varphi)$ such that
\[
\sum_{j\ge1}|D^iq_j(y)|\le C_i|y|^{-d-\sigma-i},\qquad y\neq0.
\]
\end{lemma}

\begin{proof}
Let $\gamma=d+\sigma+i$. Since $D^iq_j=2^{j\gamma}(D^iq)(2^j\cdot)$ and $q\in\cS$, we have $|D^iq_j(y)|\le N2^{j\gamma}(1+2^j|y|)^{-\gamma-1}$. The sum over $j$ with $2^j|y|\le1$ is bounded by $N\sum_{2^j\le1/|y|}2^{j\gamma}\le N|y|^{-\gamma}$. The sum over $j\ge1$ with $2^j|y|>1$ is bounded by
\[
N|y|^{-\gamma-1}\sum_{j\ge1,\,2^j>1/|y|}2^{-j}\le N|y|^{-\gamma-1}\min\{1,2|y|\}\le N|y|^{-\gamma}.
\]
The lemma is proved.
\end{proof}

For $j\ge1$, we set $L_jv(x)=\frac12\int_{\R^d}\delta_yv(x)\,q_j(y)\,dy$.

\begin{lemma}\label{lem:Ljek}
For $j\ge1$ and $k\ge0$, we have $L_jb_k=2^{k\sigma}b_k$ if $k=j$, and $L_jb_k=0$ otherwise.
\end{lemma}

\begin{proof}
We have $\frac12\delta_yb_k(x)=b_k(x)(\cos(2^ky_1)-1)$. Since $q_j$ is even and $\int q_j=0$, we get $L_jb_k=\hat q_j(2^ke_1)b_k=2^{k\sigma}\varphi(2^{k-j}e_1)b_k$. If $k\neq j$, then $|2^{k-j}e_1|\le1/2$ or $|2^{k-j}e_1|\ge2$, so that $\varphi(2^{k-j}e_1)=0$. Finally, $\varphi(e_1)=1$.
\end{proof}

\begin{lemma}\label{lem:lac}
For $p\in[1,\infty)$ and $R\ge1$, there is a constant $C=C(d,p,R)$ such that for any integers $1\le J_1\le J_2$ and real numbers $\alpha_{J_1},\ldots,\alpha_{J_2}$,
\[
\Big\|\sum_{j=J_1}^{J_2}\alpha_j\cos(2^jx_1)\Big\|_{L_p(B_R)}+\Big\|\sum_{j=J_1}^{J_2}\alpha_j\sin(2^jx_1)\Big\|_{L_p(B_R)}\le C\Big(\sum_{j=J_1}^{J_2}\alpha_j^2\Big)^{1/2}.
\]
\end{lemma}

\begin{proof}
We consider the cosine sum; the sine sum is treated in the same way. Let $F(t)=\sum_{j=J_1}^{J_2}\alpha_j\cos(2^jt)$. We first show that for any integer $n\ge0$,
\begin{equation}\label{eq:moment}
\int_0^{2\pi}F(t)^{2n}\,dt\le2\pi(2n)!\Big(\sum_{j=J_1}^{J_2}\alpha_j^2\Big)^n.
\end{equation}
We use induction on $J_2-J_1$, simultaneously for all $n$. If $J_1=J_2$, \eqref{eq:moment} is obvious. Otherwise, the shift $t\mapsto t+\pi2^{-J_1}$ changes the sign of $\cos(2^{J_1}t)$ and leaves $\cos(2^jt)$ unchanged for $j>J_1$. Thus, by periodicity, with $G(t)=\sum_{j=J_1+1}^{J_2}\alpha_j\cos(2^jt)$,
\begin{align*}
\int_0^{2\pi}F^{2n}\,dt&=\frac12\int_0^{2\pi}\Big(\big(G+\alpha_{J_1}\cos(2^{J_1}t)\big)^{2n}+\big(G-\alpha_{J_1}\cos(2^{J_1}t)\big)^{2n}\Big)dt\\
&=\sum_{i=0}^n\binom{2n}{2i}\int_0^{2\pi}G^{2i}\big(\alpha_{J_1}\cos(2^{J_1}t)\big)^{2n-2i}\,dt.
\end{align*}
By the induction hypothesis applied to $G$,
\[
\int_0^{2\pi}F^{2n}\,dt\le2\pi(2n)!\sum_{i=0}^n\frac{\alpha_{J_1}^{2(n-i)}}{(2n-2i)!}\Big(\sum_{j>J_1}\alpha_j^2\Big)^i\le2\pi(2n)!\sum_{i=0}^n\binom ni\alpha_{J_1}^{2(n-i)}\Big(\sum_{j>J_1}\alpha_j^2\Big)^i,
\]
and the right-hand side is equal to that of \eqref{eq:moment}. Given $p\in[1,\infty)$, let $n=\lceil p/2\rceil$. By H\"older's inequality and \eqref{eq:moment},
\[
\|F\|_{L_p(0,2\pi)}\le(2\pi)^{1/p}\big((2n)!\big)^{1/(2n)}\Big(\sum_{j=J_1}^{J_2}\alpha_j^2\Big)^{1/2}.
\]
Since $F$ is $2\pi$-periodic and $B_R\subset[-R,R]^d$, the lemma follows from Fubini's theorem.
\end{proof}

\begin{lemma}\label{lem:ek-lower}
For $p\in[1,\infty)$, we have $\liminf_{k\to\infty}\|b_k\|_{L_p(B_1)}>0$.
\end{lemma}

\begin{proof}
By the Riemann--Lebesgue lemma,
\[
\int_{B_1}b_k^2\,dx=\frac12|B_1|+\frac12\int_{B_1}\cos(2^{k+1}x_1)\,dx\to\frac12|B_1|.
\]
If $p\ge2$, then $\|b_k\|_{L_2(B_1)}\le|B_1|^{1/2-1/p}\|b_k\|_{L_p(B_1)}$. If $p<2$, then $|b_k|^p\ge b_k^2$ because $|b_k|\le1$.
\end{proof}

\subsection{The coefficient}
Let $\kappa\in(0,1]$ be a constant to be chosen below, and set
\[
s_j=\kappa\,\frac{\omega(2^{-j})}{\omega(1/2)},\qquad j\ge1.
\]
Then $0<s_j\le\kappa$. Since $\omega$ is nondecreasing,
\[
\omega(2^{-j})^2\log2\ge\int_{2^{-j-1}}^{2^{-j}}\frac{\omega(t)^2}{t}\,dt,
\]
and thus $\sum_js_j^2=\infty$ by \eqref{eq:notsqD}. Let $\beta_1=3$. Given $\beta_m$, let $\gamma_m\ge\beta_m$ be the smallest integer such that
\[
S_m:=\sum_{j=\beta_m}^{\gamma_m}s_j^2\ge m^2,
\]
and let $\beta_{m+1}=\gamma_m+3$. We write $I_m=\{\beta_m,\ldots,\gamma_m\}$. The blocks $I_m$ are disjoint, $\beta_m\to\infty$, and no index $\beta_m-2$ belongs to any block. We define
\[
P(x,y)=\sum_{m\ge1}b_{\beta_m-2}(x)\sum_{j\in I_m}s_j\,b_j(x)\,q_j(y),\qquad a(x,y)=c_{d,\sigma}+|y|^{d+\sigma}P(x,y).
\]

\begin{lemma}\label{lem:coef}
Let $\kappa=\min\{1,c_{d,\sigma}/(2C_0)\}$, where $C_0$ is the constant in Lemma~\ref{lem:qsum}. Then $a\in C^\infty(\R^d\times(\R^d\setminus\{0\}))$, $a$ is symmetric in $y$, and \eqref{eq:sharp-a} holds with a constant $C=C(d,\sigma,\varphi,\omega)$.
\end{lemma}

\begin{proof}
For $|y|\ge\varepsilon>0$, the series defining $P$ and all of its derivatives converge uniformly. Indeed, for multi-indices $\alpha$ and $\alpha'$ and for any $\ell>0$, we have
\[
|\partial_x^\alpha(b_{\beta_m-2}b_j)|\le N2^{j|\alpha|},\qquad|\partial_y^{\alpha'}q_j(y)|\le N_\ell2^{j(d+\sigma+|\alpha'|)}(1+2^j\varepsilon)^{-\ell}.
\] Hence $a\in C^\infty(\R^d\times(\R^d\setminus\{0\}))$, the derivatives may be taken term by term, and $a$ is symmetric in $y$ because each $q_j$ is radial.

By Lemma~\ref{lem:qsum} and $s_j\le\kappa$, we have $|D^i_yP(x,y)|\le\kappa C_i|y|^{-d-\sigma-i}$ for $i=0,1,2$. By the choice of $\kappa$, this gives $c_{d,\sigma}/2\le a\le3c_{d,\sigma}/2$, and the derivative bounds in \eqref{eq:sharp-a} follow because $a(x,y)|y|^{-d-\sigma}=c_{d,\sigma}|y|^{-d-\sigma}+P(x,y)$.

For the modulus of continuity, let $r=|x_1-x_2|$. For $j\in I_m$, we have $2^{\beta_m-2}\le2^j$, so that
\[
\big|b_{\beta_m-2}(x_1)b_j(x_1)-b_{\beta_m-2}(x_2)b_j(x_2)\big|\le\min\big\{2,(2^{\beta_m-2}+2^j)r\big\}\le2\min\{1,2^jr\}.
\]
Hence, by Lemma~\ref{lem:qsum},
\[
|a(x_1,y)-a(x_2,y)|\le2|y|^{d+\sigma}\sum_{j\ge1}s_j\min\{1,2^jr\}|q_j(y)|\le2C_0\sup_{j\ge1}s_j\min\{1,2^jr\}.
\]
We claim that $\sup_{j\ge1}s_j\min\{1,2^jr\}\le\kappa\,\omega(r)/\omega(1/2)$. If $r>1/2$, this follows from $s_j\le\kappa$ and $\omega(r)\ge\omega(1/2)$. Let $0<r\le1/2$. If $2^{-j}\le r$, then $s_j\le\kappa\,\omega(r)/\omega(1/2)$ because $\omega$ is nondecreasing. If $2^{-j}>r$, then $r<2^{-j}\le1/2$, and since $\omega(t)/t$ is nonincreasing on $(0,1/2]$,
\[
s_j2^jr=\frac{\kappa}{\omega(1/2)}\cdot\frac{\omega(2^{-j})}{2^{-j}}\cdot r\le\frac{\kappa}{\omega(1/2)}\cdot\frac{\omega(r)}{r}\cdot r=\frac{\kappa\,\omega(r)}{\omega(1/2)}.
\]
This proves the claim and the lemma.
\end{proof}

\subsection{The test functions}
We set
\[
v_m=\frac2{S_m}\sum_{j\in I_m}s_j2^{-j\sigma}b_j,\qquad u_m=2^{-(\beta_m-2)\sigma}b_{\beta_m-2}+v_m.
\]
These are finite trigonometric sums in $x_1$, so that $u_m,v_m\in C^\infty(\R^d)$ and all their derivatives are bounded. We note that the subscript $m$ here refers to the $m$-th block $I_m$, and $u_m$ should not be confused with the Littlewood--Paley pieces $u_j=\Delta_ju$ used in Sections~\ref{sec:aux} and \ref{sec:bdd}, which do not appear in this section.

\begin{lemma}\label{lem:Lum}
For any $m\ge1$, we have
\[
(-\Delta)^{\sigma/2}v_m=\frac2{S_m}\sum_{j\in I_m}s_jb_j,\qquad(-\Delta)^{\sigma/2}u_m=b_{\beta_m-2}+(-\Delta)^{\sigma/2}v_m,
\]
and
\[
\cL^au_m=-\frac2{S_m}\sum_{j\in I_m}s_jb_j+\frac{b_{\beta_m-2}}{S_m}\sum_{j\in I_m}s_j^2b_{j+1},\qquad\cL^av_m=\cL^au_m+b_{\beta_m-2}.
\]
\end{lemma}

\begin{proof}
The first two identities follow from \eqref{eq:frac-ek}. Let $v$ be a finite linear combination of the functions $b_k$. Then $|\delta_yv|\le N\min\{|y|^2,1\}$, so that by Lemma~\ref{lem:qsum} the series defining $P$ can be integrated term by term against $\delta_yv$. Since $\cL^{c_{d,\sigma}}=\partial^\sigma=-(-\Delta)^{\sigma/2}$, we obtain
\[
\cL^av=-(-\Delta)^{\sigma/2}v+\sum_{m'\ge1}b_{\beta_{m'}-2}\sum_{j\in I_{m'}}s_jb_j\,L_jv.
\]
By Lemma~\ref{lem:Ljek}, $L_jb_{\beta_m-2}=0$ for every $j$ in every block, and $L_jb_k\neq0$ only if $j=k$. Therefore,
\[
\cL^a\big(2^{-(\beta_m-2)\sigma}b_{\beta_m-2}\big)=-b_{\beta_m-2},\qquad\cL^av_m=-\frac2{S_m}\sum_{j\in I_m}s_jb_j+\frac2{S_m}b_{\beta_m-2}\sum_{j\in I_m}s_j^2b_j^2.
\]
Since $2b_j^2=1+b_{j+1}$ and $\sum_{j\in I_m}s_j^2=S_m$, we have
\[
\frac2{S_m}b_{\beta_m-2}\sum_{j\in I_m}s_j^2b_j^2=b_{\beta_m-2}+\frac{b_{\beta_m-2}}{S_m}\sum_{j\in I_m}s_j^2b_{j+1}.
\]
The lemma is proved.
\end{proof}

\begin{lemma}\label{lem:props}
Let $p\in[1,\infty)$ and $R\ge1$. Then, as $m\to\infty$,
\begin{enumerate}[label=(\roman*)]
\item $\|u_m\|_{L_\infty}+\|v_m\|_{L_\infty}\to0$;
\item $\|\cL^au_m\|_{L_p(B_R)}+\|(-\Delta)^{\sigma/2}v_m\|_{L_p(B_R)}\to0$;
\item if $\sigma>1$, then $\|Du_m\|_{L_\infty}+\|Dv_m\|_{L_\infty}\to0$; if $\sigma=1$, then $\|Dv_m\|_{L_p(B_R)}\to0$ and $\sup_m\|Du_m\|_{L_p(B_R)}<\infty$;
\item $\liminf_m\|(-\Delta)^{\sigma/2}u_m\|_{L_p(B_1)}>0$ and $\liminf_m\|\cL^av_m\|_{L_p(B_1)}>0$.
\end{enumerate}
\end{lemma}

\begin{proof}
Let $C$ be the constant in Lemma~\ref{lem:lac}. Recall that $s_j\le\kappa\le1$ and $S_m\ge m^2$.

(i) We have
\[
\|v_m\|_{L_\infty}\le\frac2{S_m}\sum_{j\ge1}2^{-j\sigma}\le\frac2{(2^\sigma-1)m^2},\qquad\|u_m\|_{L_\infty}\le2^{-(\beta_m-2)\sigma}+\|v_m\|_{L_\infty}.
\]

(ii) By Lemma~\ref{lem:Lum}, $|b_{\beta_m-2}|\le1$, Lemma~\ref{lem:lac}, and $\sum_{j\in I_m}s_j^4\le\kappa^2S_m$,
\[
\|\cL^au_m\|_{L_p(B_R)}\le\frac{2C}{S_m}\Big(\sum_{j\in I_m}s_j^2\Big)^{1/2}+\frac C{S_m}\Big(\sum_{j\in I_m}s_j^4\Big)^{1/2}\le\frac{3C}{\sqrt{S_m}},
\]
and similarly $\|(-\Delta)^{\sigma/2}v_m\|_{L_p(B_R)}\le2C/\sqrt{S_m}$.

(iii) Only $\partial_1u_m$ and $\partial_1v_m$ are nonzero, and
\[
\partial_1v_m=-\frac2{S_m}\sum_{j\in I_m}s_j2^{j(1-\sigma)}\sin(2^jx_1),\qquad\partial_1u_m=-2^{(\beta_m-2)(1-\sigma)}\sin(2^{\beta_m-2}x_1)+\partial_1v_m.
\]
If $\sigma>1$, then
\[
\|\partial_1v_m\|_{L_\infty}\le2(2^{\sigma-1}-1)^{-1}m^{-2},\qquad\|\partial_1u_m\|_{L_\infty}\le2^{(\beta_m-2)(1-\sigma)}+\|\partial_1v_m\|_{L_\infty}.
\]
If $\sigma=1$, Lemma~\ref{lem:lac} for sine sums gives $\|\partial_1v_m\|_{L_p(B_R)}\le2C/\sqrt{S_m}$ and $\|\partial_1u_m\|_{L_p(B_R)}\le|B_R|^{1/p}+2C/\sqrt{S_m}$.

(iv) By Lemma~\ref{lem:Lum} and (ii) with $R=1$,
\[
\|(-\Delta)^{\sigma/2}u_m\|_{L_p(B_1)}\ge\|b_{\beta_m-2}\|_{L_p(B_1)}-\frac{2C}{\sqrt{S_m}},\qquad\|\cL^av_m\|_{L_p(B_1)}\ge\|b_{\beta_m-2}\|_{L_p(B_1)}-\frac{3C}{\sqrt{S_m}}.
\]
Since $\beta_m\to\infty$, Lemma~\ref{lem:ek-lower} gives (iv).
\end{proof}

\subsection{Proof of Theorem~\ref{thm:sharp}}
Let $a$ be the coefficient in Lemma~\ref{lem:coef}, and let $p\in(1,\infty)$. Fix $\eta\in C_0^\infty(B_4)$ with $0\le\eta\le1$ and $\eta=1$ on $B_2$. We first note that if $g\in C^2(\R^d)$ is bounded and vanishes in $B_2$, then for $x\in B_1$, $g$ vanishes in $B_1(x)$, and thus
\begin{equation}\label{eq:far}
|\cL^ag(x)|+|\partial^\sigma g(x)|\le N\|g\|_{L_\infty},\qquad x\in B_1.
\end{equation}
Moreover, for bounded $g$, we have $\|g\|_{L_p(\mu)}\le N\|g\|_{L_\infty}$.

\emph{Proof of (c).} Suppose that $\sigma\neq1$ and the first estimate in (c) holds. We apply it to $u=u_m$. By Lemma~\ref{lem:props} (i)--(iii), the right-hand side tends to zero, which contradicts Lemma~\ref{lem:props} (iv). Now let $\sigma=1$, and set $M=\sup_m\|Du_m\|_{L_p(B_4)}<\infty$. Choose $\theta\in(0,1]$ such that $2\theta M<\liminf_m\|\partial^1u_m\|_{L_p(B_1)}$, which is possible by Lemma~\ref{lem:props} (iv). If the second estimate in (c) holds with this $\theta$ and some $N_\theta$, then by Lemma~\ref{lem:props} (i) and (ii), we obtain $\liminf_m\|\partial^1u_m\|_{L_p(B_1)}\le\theta M$, which is a contradiction.

\emph{Proof of (b).} Note that $\eta u_m\in C_0^\infty(\R^d)$ and $\|\eta u_m\|_{L_p}\le N\|u_m\|_{L_\infty}\to0$. By \eqref{eq:comm},
\[
\|\cL^a(\eta u_m)\|_{L_p}\le\|\cL^au_m\|_{L_p(B_4)}+\|u_m\,\cL^a\eta\|_{L_p}+\|I^a[\eta,u_m]\|_{L_p}.
\]
The first term tends to zero by Lemma~\ref{lem:props} (ii). By Lemma~\ref{lem:cutoff} with $R=4$, the sum of the last two terms is bounded by $N\|u_m\|_{L_p(\mu)}$ if $\sigma<1$, and by $\theta\|Du_m\|_{L_p(B_6)}+N_\theta\|u_m\|_{L_p(\mu)}$ for any $\theta\in(0,1]$ if $\sigma\ge1$. By Lemma~\ref{lem:props} (i) and (iii), and since $\theta$ is arbitrary, we conclude that $\|\cL^a(\eta u_m)\|_{L_p}\to0$. On the other hand, $\partial^\sigma(\eta u_m)=\partial^\sigma u_m-\partial^\sigma((1-\eta)u_m)$, and by \eqref{eq:far},
\[
\|\partial^\sigma(\eta u_m)\|_{L_p}\ge\|\partial^\sigma u_m\|_{L_p(B_1)}-N\|u_m\|_{L_\infty}.
\]
By Lemma~\ref{lem:props} (i) and (iv), $\liminf_m\|\partial^\sigma(\eta u_m)\|_{L_p}>0$. If \eqref{eq:global} with $w=1$ held for some $\lambda\ge0$ and $N$, then
\[
\|\partial^\sigma(\eta u_m)\|_{L_p}\le N\big(\lambda\|\eta u_m\|_{L_p}+\|\cL^a(\eta u_m)\|_{L_p}\big)\to0,
\]
which is a contradiction.

\emph{Proof of (a).} Note that $\eta v_m\in C_0^\infty(\R^d)$ and $\|\eta v_m\|_{L_p}\le N\|v_m\|_{L_\infty}\to0$. By \eqref{eq:comm} with $\partial^\sigma$ in place of $\cL^c$,
\[
\|\partial^\sigma(\eta v_m)\|_{L_p}\le\|(-\Delta)^{\sigma/2}v_m\|_{L_p(B_4)}+\|v_m\,\partial^\sigma\eta\|_{L_p}+\|I^{c_{d,\sigma}}[\eta,v_m]\|_{L_p}.
\]
By Lemma~\ref{lem:cutoff} with $R=4$ and Lemma~\ref{lem:props} (i)--(iii), the right-hand side tends to zero. Hence $\|\eta v_m\|_{H^\sigma_p}\to0$ by \eqref{eq:Hequiv}. On the other hand, $\cL^a(\eta v_m)=\cL^av_m-\cL^a((1-\eta)v_m)$, so that by \eqref{eq:far},
\[
\|\cL^a(\eta v_m)\|_{L_p}\ge\|\cL^av_m\|_{L_p(B_1)}-N\|v_m\|_{L_\infty}.
\]
By Lemma~\ref{lem:props} (i) and (iv), $\liminf_m\|\cL^a(\eta v_m)\|_{L_p}>0$. Therefore, $\cL^a$ is not bounded from $H^\sigma_p(\R^d)$ to $L_p(\R^d)$.

This completes the proof of Theorem~\ref{thm:sharp}. \qed

\begin{remark}
The proof of (a) shows that, in the example above, the failure of the a priori estimates is driven by the unboundedness of the operator $\cL^a$: the perturbation operator with kernel $P$ maps functions with small $H^\sigma_p$ norm to functions of size one. This is consistent with Theorem~\ref{thm:bdd}, since $[\omega]=\infty$.
\end{remark}

\end{document}